\documentclass{amsart}
\author{Hannah Cairns}
\usepackage{amsfonts,amssymb}\usepackage{tikz-cd}\usepackage{graphicx}
\def\angry#1.{\noindent {\bf #1.}}
\def\compatibilityfigure#1{\includegraphics{sandpile_turing_machine_fig#1.pdf}}

\title[Some halting problems for sandpiles]{Some halting problems for abelian sandpiles are undecidable in dimension three}
\newtheorem{theorem}{Theorem}
\newtheorem{lemma}{Lemma}
\newtheorem{definition}{Definition}
\newtheorem{corollary}{Corollary}

\begin{document}
\begin{abstract}
The abelian sandpile model is a simple combinatorial model for critical behaviour, a chip game on a graph with the {\it abelian property} that the order in which we make moves does not change the final outcome of the game. We prove the undecidability of three halting problems for the sandpile on the three-dimensional lattice. It is unknown whether or not the problems are decidable on a two-dimensional lattice, or on periodic planar graphs. \end{abstract}

\maketitle

\thispagestyle{empty}

\tableofcontents

\long\def\eat#1{}

\begin{section}{Introduction}

A {\sl sandpile} is a discrete process on a directed or undirected graph $G = (V, E)$. We begin with a certain number of chips on each vertex of the graph, given by a {\sl starting configuration} $S_0: V \to \mathbb Z$.

Any vertex $v$ that has at least as many chips as its out-degree $\deg^+(v)$ is said to be {\sl unstable}, and we may {\sl topple} that vertex as follows: for each outgoing edge $(v, w)$, we subtract one chip from $v$ and add one to $w$. The total effect is to take $\deg^+(v)$ chips off of $v$ and add one to each outgoing neighbour.
It makes sense to topple at a stable vertex also, and in that case the number of chips after toppling will be negative.

The process works as follows: we start with the starting configuration, topple some vertex that is unstable, and repeat. If at some point there is no unstable vertex, then the sandpile process {\sl halts}; otherwise it continues forever.

The choice of unstable vertex at each step is arbitrary, except that we insist that every unstable vertex must eventually be chosen. That is, if a vertex is unstable at some step of the process, it must be chosen at some later step, and if it is still unstable, it must be chosen again, and so on.
Notice that the only way to reduce the number of chips on a vertex is to topple it, so an unstable vertex will stay unstable until chosen.

It turns out that the arbitrary choices of vertex don't affect whether or not the process halts, or even the number of times that a particular vertex topples. This is the so-called `abelian property' of the system, which is proven in the appendix and described in the following section.

\begin{subsection}{The abelian property of the sandpile}

\label{abelianprop}Let a {\sl toppling sequence} be a finite or infinite sequence of vertices in $V$. We do not make any other restrictions: a toppling sequence may even be empty. However, not all toppling sequences are suitable for the sandpile process. We have to make sure that every vertex we topple is unstable, and we have to make sure that we eventually topple every unstable vertex.

Put chips on each vertex in $V$ according to the starting configuration $S_0$, so that the vertex $v$ starts out with $S_0(v)$ chips on it, and then topple the vertices in the toppling sequence in order, regardless of whether they are stable.
We say that the sequence is {\sl legal} if every vertex that's toppled is unstable after the previous topplings have been carried out, and {\sl complete} if every vertex that is unstable at some step is toppled at a later step. A toppling sequence is a legitimate set of choices for the sandpile process if and only if it is both legal and complete.

If $a$ is a toppling sequence, let $N_a(v) = \#\{i: a_i = v\}$ be the number of times that $v$ appears in $a$. Then we have the following three lemmas:

\begin{lemma}
Given any finite or countable graph $(V, E)$ and any sandpile configuration on it, there exists a legal complete toppling sequence for that configuration.
\end{lemma}

If $a$ is a toppling sequence, let $N_a(v) = \#\{i: a_i = v\}$ be the number of times that $v$ appears in $a$. Then:

\begin{lemma}
If $a$ is legal and $b$ is complete, $N_a(v) \le N_b(v)$ for all vertices $v$.
There is a function $N: V \to \mathbb Z_{\ge 0} \cup \{\infty\}$ so that every legal complete toppling sequence topples exactly $N(v)$ times at vertex $v$, called the odometer function.
\end{lemma}

These two lemmas are proven in Appendix \ref{abelianproperty} as Lemmas \ref{legalcomplete} and \ref{aplemma}.

\end{subsection}

\begin{subsection}{The computational power of the sandpile}
Dynamical systems typically proceed step by step, and differences in timing or speed can dramatically affect the outcome of the process.
But in sandpile dynamics, topplings can be reordered or delayed indefinitely without affecting the final state.
It's natural to ask whether a system with the abelian property has less computational power.

To ask a more concrete question, let $\mathbb Z^3$ be the simple cubic lattice graph, where every vertex has six neighbours along the cardinal axes. Is there a Turing machine that can predict the halting properties of a sandpile on that lattice? (The same methods will work on $\mathbb Z^d$ for $d > 3$, and on the slab $\mathbb Z^2 \times \{0, \ldots, n-1\}$ for $n \ge 3$.)

We consider three such problems in this paper:

\smallskip
{\it Vertex prediction.} Given a certain infinite sandpile in $\mathbb Z^3$, does the vertex at the origin ever topple?

\smallskip
{\it Global halting prediction.} Are there infinitely many topplings overall?

\smallskip
{\it Local halting prediction.} Are there infinitely many topplings at the origin?

\smallskip
First, we must decide what class of sandpiles we are going to ask about. General sandpile configurations will be too complicated, periodic or finite sandpile configurations will be too simple, and we will ultimately settle on the class of configurations that are periodic except for a finite number of added chips.
\end{subsection}

\begin{subsection}{General sandpiles are too hard}

If we allow our sandpiles to be completely general, and we encode them in a naive way, all three problems will be unsolvable because of a trivial technical issue. 
Suppose that we have a very naive encoding, where each square of the tape tells us the number of chips on some vertex.
Assume that we have a Turing machine which solves one of the problems. Run it on the sandpile where every vertex contains five chips. No vertex ever has six chips in it, so no vertex ever topples. So for this sandpile the answer to all of the problems is `no,' and our algorithm must halt at some finite time with a negative result.

However, the Turing machine only ran for a finite amount of time before halting, so it only looked at finitely many squares of its tape. That means it made its decision on the basis of the number of chips on finitely many vertices.
If we add an extra chip to one of the vertices that hasn't been seen, so that it now has six chips and is unstable, {\it every} vertex of the resulting sandpile topples infinitely often.\footnote{Let $v$ be the vertex with six chips, and let $A_i = \{w: \Vert v - w \Vert_1 = i\}$ be the set of vertices at graph distance $i$ from $v$. Topple the vertices in the following sets, one after another: $A_0$, $A_1$, $A_0$, $A_2$, $A_1$, $A_0$, $A_3$, $A_2$, $A_1$, $A_0$, \ldots. It's left to the reader to check that every toppled vertex will have six, seven, or eight chips on it: five plus the number of neighbours with smaller $\ell^1$ distance to $v$.\label{fiveorsix}} For this new sandpile, the answer to all three problems is `yes.'

Run the old machine on this new sandpile. It must halt with the same answer, `no,' because it looks at the same set of tape squares, but that answer is now wrong. So our algorithm is defective, it doesn't solve the problem after all, and we have a contradiction. No such Turing machine can exist.

More complicated encoding schemes are possible, of course, but every scheme in which a square of tape gives us information about only finitely many vertices will fail for the same reason. Also, we are going to see later that the halting problem for Turing machines is reducible to a special case of this problem, so the combination of the encoder and solver has to be able to solve the halting problem.

The question is not really fair. A terminating Turing machine isn't suited to solving problems with an infinite amount of input data, because it can only read a finite amount of that data before making its decision. How can we modify this problem to make the input data finite and give the Turing machine a chance?

\end{subsection}

\begin{subsection}{Periodic sandpiles are too easy}\label{periodicsandpile}
One natural thing to do is to restrict to periodic sandpile configurations, but that special case is too easy: all the problems have trivial solutions, if perhaps not computationally fast ones.

A sandpile configuration on $\mathbb Z^3$ is {\sl periodic} if there are integers $n_x, n_y, n_z \ge 1$ so that two vertices $(x,y,z)$ and $(x',y',z') \in \mathbb Z^3$ have the same number of chips whenever~\[x \equiv x' \bmod n_x, \ y \equiv y' \bmod n_y, \ z \equiv z' \bmod n_z.\]
If we know the three periods, and the number of chips on every vertex in the origin cube $C_0 = \{(x, y, z) \in \mathbb Z^3 \mid 0 \le x < n_x,\, 0 \le y < n_y,\, 0 \le z < n_z\}$, that determines the number of chips on every other vertex in $\mathbb Z^3$.\footnote{Strictly speaking, we should call it the ``origin rectangular cuboid'' because the side lengths may be unequal, but for the sake of brevity we'll refer to it and other things like it as cubes.} That's a finite amount of information, so it's possible to encode it in finitely many tape squares.

We ask the three questions for periodic sandpiles:

\smallskip
{\it Periodic vertex prediction.} Given a certain infinite {\it periodic} sandpile in $\mathbb Z^3$, does the vertex at the origin ever topple?

\smallskip
{\it Periodic global halting prediction.} Given a certain infinite periodic sandpile in $\mathbb Z^3$, are there infinitely many topplings overall?

\smallskip
{\it Periodic local halting prediction.} Given a certain infinite periodic sandpile in $\mathbb Z^3$, are there infinitely many topplings at the origin?

\smallskip
These problems are easy to solve once we have the following lemma.

\begin{lemma}If $N(v)$ is the odometer function of a periodic sandpile configuration $S_0$, and $N_T(v)$ is the odometer function of the corresponding configuration on the finite torus graph $T$, then $N(v) = N_T([v])$ for every vertex $v \in \mathbb Z^3$.
\end{lemma}

This is proven as Lemma~\ref{sameastorus} in Appendix~\ref{abelianproperty}.
Once we know this, we can solve all three problems by finding the odometer function $N_T: T \to \mathbb Z_{\ge 0} \cup \{\infty\}$ for the torus graph. This can be done in finite time since there are only finitely many possible states for the sandpile\footnote{We can only topple unstable vertices, so a vertex can't drop below zero chips, and if it starts with a negative number of chips, it won't go below that number. Toppling doesn't alter the total number of chips. Therefore the number of chips on each vertex is bounded below and above.}. We omit the details.

Once we have the odometer function on the torus, we can read off the answers to the halting questions.
\smallskip

\begin{itemize}
\item Vertex prediction: the origin topples if and only if $N_T(\vec 0) > 0$.
\item Global halting: there are infinitely many topplings if and only if $N_T \not\equiv 0$.
\item Local halting: the origin topples exactly $N_T(\vec 0)$ times.
\end{itemize}

All three problems for periodic sandpiles are solvable by a Turing machine, and the solution is straightforward, so the class of periodic sandpiles is too small to produce interesting questions about undecidability. The problem is that the periodic sandpile is translation symmetric, which makes it equivalent to a sandpile on a finite graph.

We will try to solve these problems for a periodic sandpile configuration with finitely many chips added to break the symmetry. First, though, an aside about the computational complexity of sandpile stabilization on finite graphs.
\end{subsection}
\begin{subsection}{Aside: The odometer function can be found in polynomial time if the starting configuration is nonnegative}

\smallskip
A random-access machine can find the odometer function for a finite graph with $n$ vertices, $m$ edges, and diameter $d$ in~$O(n^2md)$ time and $O(n+m)$ space, as long as each vertex starts with $\ge0$ chips.
For the torus above, that means $O(n_x^3n_y^3n_z^3(n_x+n_y+n_z))$ time and $O(n_xn_yn_z)$ space.

\begin{theorem}If the starting configuration is nonnegative, a sandpile on a finite simple undirected graph with $n$ vertices halts after $\le2nmd$ topplings or never halts.\label{sandpileatmost}\end{theorem}

\begin{proof}This is Theorem 2 in Tardos~\cite{gabortardos}.
\end{proof}

Break the graph up into graph components $G_1, \ldots, G_C$ with $n_1, \ldots, n_C$ vertices, $m_1, \ldots, m_C$ edges. That takes $O(n+m)$ time with Tarjan's algorithm. Now we want to determine whether each component runs forever or eventually halts, and if it halts, we want to know how many times each vertex topples.\footnote{If a sandpile on a connected finite graph doesn't halt, then every vertex topples infinitely often. Otherwise there would be an edge between a vertex that topples infinitely many times and a vertex that topples finitely many times, and the second vertex would receive infinitely many chips after its last toppling, which is impossible by completeness.}

Count the number of chips in each component $G_i$. If there are more than $2m_i - n_i$, some vertex will always be unstable, so $G_i$ can never stabilize.

Otherwise, run the sandpile process it until it halts or until we have toppled $n_im_id$ times, keeping track of how many times each vertex topples.
{} There are at most $n_i-1$ edges incident to a vertex, so each toppling takes $O(n_i)$ time.
\ Let the final count of topplings be $M(v)$ for $v \in V$.

The graph components that don't halt continue to topple forever, so each component contains at least one vertex that topples infinitely many times.

A vertex that topples infinitely often will put an infinite number of chips on its neighbours, which by completeness means that they topple infinitely often also. Each graph component is connected, so all vertices in a non-halting component topple infinitely many times, and $N(v) = \infty$ in those components.
%

In components that halt, the odometer function is obviously equal to $M(v)$.
Therefore, the odometer function is $$N(v) = \begin{cases}M(v)&\text{if the graph component with $v$ in it has halted}\\\infty&\text{otherwise}.\end{cases}$$

Each toppling will take $O(n)$ time because there may be as many as $n-1$ edges incident to a vertex, so the algorithm takes $O(n + m + n_1^2m_1d + \cdots + n_C^2m_Cd) = O(n^2md)$ time. We have to record the graph, the current number of chips on each vertex, and the odometer function $M: T \to \mathbb Z_{\ge 0}$, which takes $O(m+n)$ space.

\end{subsection}

\begin{subsection}{We will work with periodic+finite sandpiles}
\label{haltingmentionedsection}
We have seen that an arbitrary sandpile is too hard to predict, and a periodic sandpile is too easy. A sandpile with a finite number of chips is also too easy: any sandpile with finitely many chips on~$\mathbb Z^d$ will
always halt (a special case of Lemma \ref{dminustwoalwayshalt} in the appendix).

Start with a periodic sandpile and add a finite number of chips. This will break the translation symmetry, so we can't determine its behaviour by looking at a finite torus. And both the periodic sandpile and the locations of the finite number of chips can be encoded in finitely many tape squares, so the problem shouldn't be too hard.

We will call a periodic sandpile which has had finitely many chips added to it a {\it periodic+finite sandpile}, or p+f sandpile.
And we ask the same questions:

\smallskip \label{problemdefinition}
{\it P+f vertex prediction.} Given a certain infinite periodic sandpile with finitely many chips added to it, does the vertex at the origin ever topple?

\smallskip
{\it P+f global halting prediction.} Given a p+f sandpile, are there infinitely many topplings overall?

\smallskip
{\it P+f local halting prediction.} Given a p+f sandpile, are there infinitely many topplings at the origin?

\smallskip
All three of these are again undecidable, but the reasons are more interesting. The rest of the paper is dedicated to the proofs.

\end{subsection}

\begin{subsection}{Outline of the proofs}

We will begin by recapitulating old material. In the 1990s, Moore and Nilsson~\cite{moorenilsson} came up with a way to embed wires, {\sc and} gates, and {\sc or} gates in sandpiles of dimension two or more. In section~\ref{howgates}, we will see what these are and how they work, and in sections~\ref{buildingturing}--\ref{vertexhalting} we will build a sandpile circuit that simulates a Turing machine as a warm-up and use it to prove that the vertex prediction problem is undecidable.

The novelty of this paper is in the adaptations of the basic idea to show that global halting and local halting are also undecidable.
In sections \ref{lazy} and \ref{evenlazier}, we will define a `lazy automaton,' built in such a way that only finitely many cells will activate if the Turing machine it is simulating halts, and use it to show that global halting is undecidable. Finally, in section \ref{localhalt}, we will build an enormous `bomb' which, if set off, will cause \emph{every} vertex in the whole sandpile to topple infinitely many times. We will connect that to our Turing machine to show that local halting is undecidable.
\end{subsection}

\end{section}

\begin{subsection}{History of halting problems on sandpiles}

We were not able to find any other work that considers the undecidability of sandpile halting problems on an infinite graph, but the computational complexity of the halting question on a finite graph without a sink has been widely considered.\footnote{Every connected undirected graph or strongly connected directed graph with a sink halts, so the halting problem for such graphs is trivial. There are other interesting problems, like computing the final state, or finding the odometer function, but they are outside the scope of this paragraph. For example, in the paper \cite{moorenilsson} mentioned above, Moore and Nilsson proved that a finite $d$-dimensional lattice of side length $n$ with a sink connected to every boundary vertex, where no vertex starts with more than $4d$ chips, will halt after no more than $O(n^{d+2})$ topplings.}

The first step was the previously mentioned result from the paper by Tardos in 1988 \cite{gabortardos}, giving an upper bound of $2nmd$ topplings for finite simple undirected graphs that halt. Bj\"orner, Lovasz, and Shor \cite{bjornerlovaszshor} gave a different polynomial bound on the number of topplings in 1991 based on the smallest eigenvalue of the adjacency matrix. In the same year, Eriksson \cite{nopolynomial} gave an example of a {\it directed} multigraph which halts after an exponential number of topplings.

The first two papers give us bounds on the complexity of halting prediction on simple finite undirected graphs, because we can always just simulate the sandpile until it halts or exceeds the bound, taking $O(\deg^+(v)) = O(n)$ time per toppling. The third paper tells us that that won't work well for directed multigraphs. Bj\"orner and Lovasz \cite{bjornerlovasz} gave a more sophisticated algorithm for directed graphs in 1992, although that algorithm also takes exponential time in the worst case.

That was the state of the art for some time. A recent paper by Farrell and Levine~\cite{farrelllevine} has shown that the sandpile halting problem is NP-complete for general finite directed multigraphs, and, on the other hand, that it can be solved in linear time for the special case of `coEulerian' directed graphs.

Of course, there may be a faster approach than direct simulation. Moore and Nilsson~\cite{moorenilsson} showed that sandpile prediction on a line is in ${\bf NC}^2$, and gave an $O(n)$ serial algorithm where $n$ is the length of the line\footnote{We consider arithmetic operations on $O(\log n)$-bit numbers to be $O(1)$ work instead of $O(\log n)$, so this differs from the complexities given in the paper, which are ${\bf NC}^3$ and $O(n \log n)$.}. Ramachandran and Schild \cite{ramachandranschild} are able to solve the prediction problem for sandpiles for the more general case of undirected trees in~$O(n \log^5 n)$ time, which can be improved to $O(n \log^2 n)$ \cite{mygithub}.




\end{subsection}

\begin{section}{Sandpiles give you fuse circuits}\label{howgates}

\begin{subsection}{Fuse circuits}\label{oneshotsection}

Following~\cite{moorenilsson}, we say that a {\sl fuse circuit} is a logical circuit built from {\sc and} and {\sc or} gates with multiple input wires and multiple output wires. No {\sc not} gates are allowed.
It's similar to a monotone circuit, but the graph of a fuse circuit is allowed to have cycles, and
to prevent oscillations or indeterminacy, we use wires with `fuses': if a wire ever turns on at any point, the fuse breaks and it stays on forever. So the state of the system is monotone and converges to a limit.

We formalize the dynamics, which are essentially asynchronous. Let $W$ be a finite or countable set of wires. Let~$\mathcal F(W)$ be the set of finite subsets of $W$.

An {\sc and} gate $a\in{}A$ is a tuple $a = (a_{\text{input}}, a_{\text{output}}) \in \mathcal{F}(W) \times 2^W$, where $a_{\text{input}}$ is a finite set of input wires and $a_{\text{output}}$ is an arbitrary set of output wires.
An {\sc or} gate is a tuple $o = (o_{\text{input}}, o_{\text{output}}) \in 2^W \times 2^W$, where both sets may be infinite.
Choose a set of {\sc and} gates $A \subseteq \mathcal F(W) \times 2^W$ and {\sc or} gates $O \subseteq 2^W \times 2^W$. This is the circuit.

If $S \subseteq W$ is the set of wires that are on at a certain step, then we \emph{may} activate any {\sc and} gate with $a_{\text{input}} \subseteq S$, and any {\sc or} gate so that $o_{\text{input}} \cap S$ is not empty. When a gate activates, we \emph{may} turn on any or all of its output wires. Or we may not. The only restriction, again, is that any wire which may be turned on is eventually turned on at some step.
Let~$\mathcal{D}(S) \subseteq W$ be the set of wires that are on at the current step, or that may be turned on at the next step:
\begin{align*}\mathcal{D}(S) = {}&S \ \cup\  \{w: \exists a \in A: a_{\rm input} \subseteq S\ \text{and}\ w \in a_{\rm output}\}\\&\qquad\qquad \cup\ \{w: \exists o\in O: o_{\rm input} \cap S \ne \varnothing\ \text{and}\ w \in o_{\rm output}\}.\end{align*}

The function $\mathcal{D}: 2^W \to 2^W$ is monotone and increasing.
Choose an initial condition $S_0 \subseteq 2^W$, where $w \in S_0$ if and only if $w$ is turned on at the start. The final state of the fuse circuit is the smallest fixed point of $\mathcal{D}$ above $S_0$, which we call $S_\infty$.

It's not difficult to see that $S_\infty$ is the union of all sets $\mathcal{D}^{(i)}(S_0)$. We can think of the circuit as a dynamical system where the wires start in the state $S_0$ and then, at each step, every gate that may be activated is activated: $S_1 = \mathcal{D}(S_0)$, $S_2 = \mathcal{D}(\mathcal{D}(S_0))$, etc. The final state is the fixed point $S_\infty$.

However, that can be misleading, because a fuse circuit is asynchronous, in the sense that timing is irrelevant. We can delay the activation of any gate until a later step without affecting the final state. As long as the current set of active wires $S$ is somewhere between $S_0$ and $S_{\infty}$, we'll reach the same fixed point $S_\infty$.

Sandpiles have the same sort of asynchronous behaviour, where the exact timing of events has no effect on the outcome. We'll see in this section that it is possible to build an arbitrary fuse circuit in a sandpile for $d \ge 3$.




\end{subsection}

\begin{subsection}{Sandpile wires}
We describe a method for building fuse circuits in a sandpile on $\mathbb Z^3$ due to Moore and Nilsson~\cite{moorenilsson}. 
We start with the wire.

A sandpile `wire' is a connected path of vertices in~$\mathbb Z^3$ with $2d-1=5$ chips on each vertex. If we add a chip to one of the vertices, it will become unstable, topple, and put chips on its neighbours, eventually toppling the whole path as in Figure~\ref{toppling}.

\begin{figure}[h]
{\arraycolsep=1pt 
\begin{center}\small
\begin{tikzcd}[ampersand replacement=\&, row sep=0em]
\parbox{4cm}{\centering A chip is added to the underlined vertex:}\&{\begin{array}{cccccccccc}\phantom{0}&5&5&5&5&5&5&\cdots\phantom{0}\\\end{array}}\&{\begin{array}{cccccccccc}\cdots&5&5&5&5&5&\underline{6}&\phantom{0}\\\end{array}}\\ \parbox{4cm}{\centering The rightmost vertex topples, putting one chip on each  neighbour:}\&{\begin{array}{cccccccccc}\phantom{0}&5&5&5&5&5&5&\cdots\phantom{0}\\\end{array}}\&{\begin{array}{cccccccccc}&&&&&&1\\\cdots&5&5&5&5&6&\underline{0}&1\\&&&&&&1\end{array}}\\
\parbox{4cm}{\centering The next-rightmost vertex topples, putting one chip on its own neighbours:}\&{\begin{array}{cccccccccc}\phantom{0}&5&5&5&5&5&5&\cdots\phantom{0}\\\end{array}}\&{\begin{array}{cccccccccc}&&&&&1&1\\\cdots&5&5&5&6&0&\underline{1}&1\\&&&&&1&1\end{array}}\\
\parbox{4cm}{\centering Eventually every vertex in the wire topples, and the final state is:}\&{\begin{array}{cccccccccc}&1&1&1&1&1&1\\1&0&1&1&1&1&1&\cdots\\&1&1&1&1&1&1\\\end{array}}\&
{\begin{array}{cccccccc}
&1&1&1&1&1&1\\
\cdots&1&1&1&1&1&\underline{1}&1\\
&1&1&1&1&1&1\end{array}}\\[-1em]
\end{tikzcd}\end{center}
}
\caption{A wire is set off. This page is only two-dimensional, so we are only showing a single horizontal slice of $\mathbb Z^3$, but there are two more lines of 1s above and below the page.}\label{toppling}
\end{figure}



\vskip .6em

The toppling of one wire will not affect another wire if they aren't adjacent. Also, we can add bends and three-way branches.
We will later have reason to be careful about how many chips we add to vertices that aren't in the wire, so let's observe now that, as long as we leave enough space between bends and branches in a wire, it will add at most two chips to any vertex that isn't part of a wire or gate.

A wire in a fuse circuit can only turn on one time, and after that it stays on forever. A sandpile wire can only topple once. We say that a wire in the fuse circuit turns on when its corresponding wire in the sandpile topples.

In a fuse circuit, the {\sc and} gate waits for all its input wires to turn on, and then turns on the output wire. So in the sandpile, a two-input, one-output {\sc and} gate should wait for both its input wires to topple and then topple its output wire. Similarly, a two-input, one-output {\sc or} gate in the sandpile should wait for one input wire to topple and then topple its output wire. We will build such gates in the next two sections.

\end{subsection}

\begin{subsection}{Wait gates and {\sc and} and {\sc or} gates} We begin with the so-called {\sc wait-2} and {\sc wait-1} gates shown in Figure \ref{2gat}.

\begin{figure}[h]\small
\begin{center}
\begin{tabular}{c|c}
\arraycolsep=1pt
$\begin{array}[b]{ccccccccc}
\cdots&5&5&5&4&5&5&5&\cdots\\
&&&&5&&&\\
&&&&5&&&
\end{array}$
&
\arraycolsep=1pt
$\begin{array}[b]{ccccccccc}
\cdots&5&5&5&5&5&5&5&\cdots\\
&&&&5&&&\\
&&&&5&&&
\end{array}$
\end{tabular}
\end{center}
\caption{{\sc wait-2} and {\sc wait-1} gates.}\label{2gat}
\end{figure}

The {\sc wait-2} gate waits for two wires to topple, and then topples the third wire, whichever one that is. The {\sc wait-1} gate waits for one wire to topple, and then topples the other two.
This isn't what we want, of course.
The logical {\sc and} gate waits for its {\it input} wires to topple and then topples the {\it output} wire, without affecting the other input wire. Similarly, the {\sc or} gate waits for one input wire and then topples the output wire.
We need some way to keep the gate from affecting the input wires.

This is accomplished by the {\sc diode} below. 
\begin{figure}[h]\small
\begin{tabular}{ccc}
\arraycolsep=1pt
$\begin{array}[b]{cccccccc}
&&&5&5&&&\\
\cdots&5&5&5&\underline{4}&5&5&\cdots\\[0em]&\\[-0.7em]
\end{array}$
&$\begin{array}{c}\\[-4em]\parbox{5cm}{\raggedleft If the left wire topples,\\the right wire topples}\  \ \ \ \ \ \ {\longrightarrow}\end{array}$&
\arraycolsep=1pt
$\begin{array}[b]{cccccccc}
&&&1&1&&&\\
&1&2&1&1&2&1&\\
\cdots&1&1&2&\underline{1}&1&1&\cdots\\
&1&1&1&1&1&1&
\end{array}$\\[1em]
&$\begin{array}{c}\\[-4em]\parbox{5cm}{\raggedleft but if the right wire topples,\\the left one doesn't}\  \ \ \ \ \ \ {\longrightarrow}\end{array}$&
\arraycolsep=1pt
$\begin{array}[b]{cccccccc}
&&&&&&&\\
&&&5&5&1&1&\\
\cdots&5&5&5&\underline{5}&0&1&\cdots\\
&&&&&1&1&
\end{array}$
\end{tabular}\label{2dio}
\caption{Diode}
\end{figure}

If the left part of the diode topples, the vertex with 4 chips will receive two more chips and will topple. That will topple the right part of the diode.
On the other hand, if the right part of the diode topples, the vertex with~4 chips receives only one chip and will not topple, leaving the left side of the diode untouched.
So signals pass from left to right through the diode, but not the other way around.

If we attach two diodes to the input wires of a {\sc wait-2} gate, as in Figure~\ref{anddiodes}, it will turn it into an {\sc and} gate. 
\begin{figure}[h]

\begin{center}
\vskip .6em \small
\arraycolsep=1pt
$\begin{array}[b]{ccccccccccccccccccc}
&&5&5&&&&5&5\\
\cdots&5&5&4&5&\underline{4}&5&4&5&5&\cdots\\
&&&&&5&&&\\
&&&&&5\\
&&&&&\vdots
\end{array}$
\end{center}
\caption{{\sc Wait-2} gate + diodes = {\sc and} gate}\label{anddiodes}\end{figure}
If the left or right input wire topples, then one chip will be added to the underlined~4 vertex.
If both input wires topple, the vertex with 4 chips in the middle will get two chips and topple the bottom output wire.

On the other hand, if two of the other wires topple --- say, if the left and bottom wires topple --- the central area will topple, but the resulting chain reaction will be absorbed by the business end of the right-hand diode and will not set off the right wire. This is shown in Figure \ref{leftbottomand}.

\begin{figure}[h]

\begin{center}
\compatibilityfigure{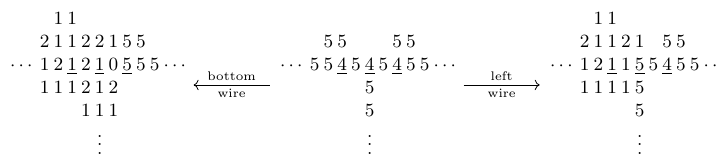}
\end{center}
\caption{An {\sc and} gate after the left and bottom wire topple}\label{leftbottomand}
\end{figure}

So this gate operates like a two-input {\sc and} gate. We could put a diode on the output wire too, but that is unnecessary, since it would only prevent the output wire from setting off the input wires, and the input diodes do that already.

If we attach two diodes to a {\sc wait-1} gate, it becomes an {\sc or} gate. See Figure~\ref{ordiodes}.
\begin{figure}[h]\begin{center}

\compatibilityfigure{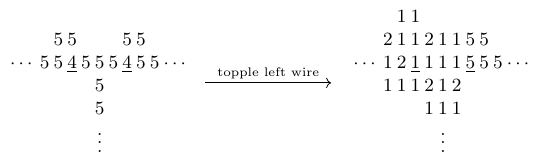}\end{center}
\caption{{\sc Wait-1} gate + diodes = {\sc or} gate}\label{ordiodes}\end{figure}

The right side of the figure shows what happens after we topple the left wire. The two~4 vertices are underlined on both sides.

The toppling traveled through the left diode in the correct direction and set off the bottom wire. The 4 square in the diode on the right received one chip, but it was in the wrong direction, so that wasn't enough to make it topple.

Having built two-input gates with one output, we can  construct {\sc and} and {\sc or} gates with arbitrarily many inputs by chaining together two-input gates as usual, and we can connect the one output wire to multiple output wires via diodes. 
(The diodes prevent output wires from toppling each other.)

\end{subsection}

\begin{subsection}{Goles and Margenstern's construction}
We note that Bitar, Goles, and Margenstern have also constructed a sandpile Turing machine in \cite{bitargoles, golesmargenstern}.
Their design works with a sandpile process that has time steps, where each unstable site at a certain time topples simultaneously, and they encode bits in the relative timing of two signals traveling through two parallel wires.

Unfortunately, that is unsuitable for the problems we want to solve, because our problems depend only on which vertices topple and how many times they topple. Their answers cannot depend on the timing of the process. In fact, by the abelian property, we can delay any toppling as long as we want without changing that data, so any bit can be flipped by delaying the leading signal until the trailing signal passes it, and the answers to our problems do not depend on the bits traveling through the computer.
That makes that design useless for these particular problems.
\end{subsection}

\begin{subsection}{Three vs.~two dimensions and a conjecture based on results from complexity theory}
\label{twod}
We are in three dimensions, so there is no planarity obstruction: given sufficient space, we can connect gates together in any arrangement we like.
This makes it possible to construct a sandpile corresponding to any fuse circuit.

In two dimensions, there is no obvious way to run sandpile wires through each other without interference, and this seems to impose a restriction of planarity.

Gajardo and Goles \cite{gajardo} prove a version of this planarity restriction for sandpiles in $\mathbb Z^2$. They show that one cannot construct a `cross-over,' a sandpile configuration that allows wire signals to cross each other without interfering. We give a variation of their argument as Lemma~\ref{crossover} in the appendix.
They also show that a crossover exists if the vertex neighbourhood is expanded to a 5x5 cross, \lower0.1em\hbox{\includegraphics[width=1em]{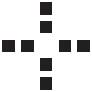}}.

Planarity is known to change the complexity behaviour of monotone circuits, assuming that ${\bf NC} \ne {\bf P}$. The evaluation of finite monotone circuits is ${\bf P}$-complete, while
fast parallel algorithms are known for the evaluation of finite {\it planar} monotone circuits.
Yang \cite{yang} shows that the evaluation of planar monotone circuits is in ${\bf NC}^3$, and on the practical side, there is an algorithm of Ramachandran and Yang \cite{ramayang} which runs in $O(\log^2 n)$ time on a linear number of parallel processors.
So, if ${\bf NC} \ne {\bf P}$, there is a qualitative difference between three and two dimensions.

On the basis of that qualitative difference for monotone circuits, we conjecture that there is also a difference in kind for infinite sandpiles:

\smallskip

{\it Conjecture.} The halting problem for periodic {\it planar} sandpiles with finitely many chips added or removed is decidable.

\smallskip
A weaker conjecture is that the halting problems for periodic+finite sandpiles in~$\mathbb Z^3$ cannot be reduced to the halting problems for periodic+finite sandpiles in $\mathbb Z^2$, so that the latter system is strictly computationally weaker, even if not decidable.

Following Bond and Levine \cite{blev}, we say that a {\sl unary network} is a message-passing network with only one type of message, where each node is a finite state machine. Such a network is necessarily abelian. The sandpile is a special case of a unary network: on $\mathbb Z^2$, each node counts the messages it gets until it hits four (or fewer if it starts with chips on it), sends one message to each of its neighbours, and repeats.

A unary network on a planar graph seems to be subject to the same kind of planarity restriction as a sandpile, so we propose the more general conjecture:

\smallskip
{\sc Conjecture.} The halting problem for a periodic unary network on $\mathbb Z^2$ with finitely many messages added is decidable.

\smallskip

That is, we conjecture that the planarity of the graph restricts the system enough that its behaviour becomes predictable.

\begin{subsubsection}{What about more complicated planar networks?\nopunct}

What if we let the message-passing network use two different types of messages instead of just one? Can we expect the halting problem for that class of networks to be decidable?

No. All our gates can be constructed in two dimensions, so the only obstacle to building an arbitrary periodic circuit is getting wires to cross each other.
With two messages, we can build wires that can cross, and get around the planarity. This is true even when the network is restricted to the class of {\sl abelian networks} (as in \cite{blev}).

Start with the two-dimensional sandpile on $\mathbb Z^2$ and replace all the vertices in the set $\{(x, y) \in \mathbb Z^2: x \bmod 8 \equiv y \bmod 8 \equiv 0\}$ with {\sl two-counter cross-over vertices}. These are vertices with two counters, one for chips that come in vertically and one for chips that come in horizontally.
When the vertical counter reaches two, it sends one chip north and one chip south. Similarly, when the horizontal counter reaches two, it sends one chip west and one chip east.

The remaining vertices behave as normal: they have one counter, and when it reaches four, they send one chip in each direction.
This is a simple system with two types of messages: vertical and horizontal chips. The message processors are abelian in the sense of Bond and Levine, so the whole system has the global abelian property by Lemma 4.4 in \cite{blev}.
But now we can build a cross-over gate as in Figure~\ref{crossove}.

\begin{figure}[h]
\begin{center}\arraycolsep=1pt
\footnotesize\parbox[c]{6em}{$\begin{array}[b]{ccccccc}
&&\vdots\\
&&3\\
&&3\\
\cdots&3&\,{}^{1}\!\!/\!{}_{1}&3&\cdots\\
&&3\\
&&3\\
&&\vdots
\end{array}$}
\parbox[c]{1.9em}{\quad}
\parbox[c]{6em}{$\begin{array}[b]{ccccccc}
&&\vdots\\
&&3\\
&&3\\
\cdots&4&\,{}^{1}\!\!/\!{}_{1}&3&\cdots\\
&&3\\
&&3\\
&&\vdots
\end{array}$}
\parbox[c]{1.9em}{$\longrightarrow$}
\parbox[c]{6em}{$\begin{array}[b]{ccccccc}
&&\vdots\\
&&3\\
&1&3&&\\
\cdots&0&\,{}^{2}\!\!/\!{}_{1}&3&\cdots\\
&1&3&&\\
&&3\\
&&\vdots
\end{array}$}
\parbox[c]{1.9em}{$\longrightarrow$}
\parbox[c]{6em}{$\begin{array}[b]{ccccccc}
&&\vdots\\
&&3\\
&1&3&&\\
\cdots&1&\,{}^{0}\!\!/\!{}_{1}&4&\cdots\\
&1&3&&\\
&&3\\
&&\vdots
\end{array}$}
\parbox[c]{1.9em}{$\longrightarrow$}
\parbox[c]{6em}{$\begin{array}[b]{ccccccc}
&&\vdots\\
&&3\\
&1&3&1&\\
\cdots&1&\,{}^{1}\!\!/\!{}_{1}&0&\cdots\\
&1&3&1&\\
&&3\\
&&\vdots
\end{array}$}
\end{center}
\caption{A crossover gate. A chip on the west wire topples the east wire without affecting the north or south wires. Here $\,{}^{1}\!\!/\!{}_{1}$ indicates a crossover vertex that starts with one vertical chip and one horizontal chip. }\label{crossove}
\end{figure}
This allows us to run wires wherever we like without a planarity obstruction. The gates described in section~\ref{howgates} all work fine in two dimensions, although the numbers have to be changed: 5 and 4 become 3 and 2. So our Turing machine constructions in later sections carry over to this situation\footnote{
The proof for local halting doesn't carry over, since the `bomb' gadget doesn't work in two dimensions: it is either not sensitive enough to explode, or so sensitive that a wire with a bend in it will set it off. Instead, we notice that the non-crossover sites are not using the second message and give them some explosive behaviour when they receive that second message: say, when they get it, they send one copy of the second message in all directions.}, and
the halting problems are undecidable for periodic two-message abelian networks on $\mathbb Z^2$, let alone more complicated ones.

Here is another way to reach this conclusion. In our construction, the graph~$\mathbb Z^3$ is broken up into a lattice of large $N \times N \times N$ cubes, and the only connections are horizontal: no cube is ever connected to a cube above or below. The sandpile will behave the same way if we delete all the vertices outside $\mathbb Z^2 \times \{0, \ldots, N\}$.

Think of each stack $\{x\} \times \{y\} \times [0, N]$ as a single node in a message processor network with $4N$ types of messages on the graph $\mathbb Z^2$.
The abelian property is not injured by this grouping procedure. The resulting network is abelian, and it's periodic except that the starting state is altered at finitely many vertices.

This network behaves identically to the Turing machine sandpile that we construct later. An algorithm to predict halting on this class of networks would allow us to predict the halting of the sandpile, so there is no such algorithm, and
again we see that vertex and global halting are undecidable for general abelian networks on $\mathbb Z^2$.

By the way, here is a word of caution. It's an interesting exercise to try to set up a Turing machine with one abelian message processor\footnote{As in Bond and Levine~\cite{blev}, an abelian message processor is a message processor with the property that reordering input messages doesn't affect the final state or the messages that are sent out, except that the outputs are allowed to be in a different order.} per tape cell per unit of time and all the processors identical, but most of the work in constructing the state machine goes into answering questions like ``What should I do if someone sends me many contradictory messages about the state of the tape cells at the previous time?'' even though the processor won't get into that situation under normal operation.

We use sandpile circuits to avoid having to think about these strange situations. They're automatically abelian, so we don't need to do any work at all to keep things consistent. (But it's worthwhile to keep the question above in mind as we build our sandpile circuits. Later we'll use pairs of wires to represent single Boolean bits. Once we get there, here's an exercise: what happens if both wires in each bit are activated?)

\end{subsubsection}

\end{subsection}

\begin{subsection}{The course of the rest of the paper}

\def\smaller{\hskip -3em}
\begin{figure}[t]\small
\centering
\begin{tikzcd}
\smaller\text{The sandpile operates like}\smaller\ar{d}{\text{\ sec.~\ref{oneshotsection}}}\\[1em]\smaller\parbox{3cm}{\centering a fuse circuit, which simulates}\ar{dd}{\text{\ sec.~\ref{buildingturing}}}\ar[rd, start anchor={south east}, end anchor={[xshift=3em]north west}, "\text{\ sec. \ref{lazy}}"]\smaller\\[2em]&
\parbox{4cm}{\centering a lazy cellular automaton which simulates}\ar[ld, start anchor={[xshift=3em]south west}, end anchor={north east}, "\text{\ sec. \ref{evenlazier}}"]
\\[2em]
\text{a Turing machine}
\end{tikzcd}
\caption{The layers of the paper.}\label{twoway}
\end{figure}

Having described Moore and Nilsson's construction of sandpile circuits and made some grand conjectures based on the work of Yang, we now get around to the actual content of this paper, three relatively prosaic proofs that the three problems earlier are undecidable.

The structure of the paper is summarized in Figure~\ref{twoway}.
We have already seen how to construct fuse circuits in sandpiles in section~\ref{oneshotsection}, which covers the top-left arrow. For vertex and local halting, we need the lower arrow on the left: we simulate a Turing machine directly with a fuse circuit, explained in the next section, section~\ref{buildingturing}.

The undecidability of vertex halting is covered in subsection~\ref{vertexhalting}.
For global halting, we introduce the idea of a {\sl lazy automaton} in section~\ref{lazy} and construct a Turing machine out of one in section~\ref{evenlazier}, with the property that the automaton does a finite amount of work if and only if the Turing machine halts eventually. Local halting is covered in section~\ref{localhalt}. We build a device that, when toppled, causes every vertex in $\mathbb Z^3$ to topple infinitely often, and show that fuse circuits can coexist with it and set it off.

\end{subsection}
\end{section}

\begin{section}{Building Turing machines from fuse circuits}\label{buildingturing}

For all three problems, we show undecidability by reducing them to the halting problem for Turing machines using some sort of special gadget which varies with the problem.

The first step in the vertex and local halting problems is to build a fuse circuit that simulates a Turing machine. There is no novelty in this part: a computer scientist could do this blindfolded, and we encourage the reader to skip to Section~\ref{lazy} if all of this seems obvious.
We start by recalling what a Turing machine is.

\begin{subsection}{Turing machines}

A Turing machine is a six-tuple $(Q, q_0, F, \Gamma, \beta, \delta)$.

\smallskip
\begin{itemize}
\item It has a set of states $Q$, one of which is the starting state $q_0$. Some subset of these states are final states $F \subseteq Q$. If the machine enters any of these states, it halts, and if it never enters any final state, it continues forever.

\item It has an infinite one-dimensional tape with one square for every $n \in \mathbb Z$,
and each square has one letter written on it from the tape
alphabet $\Gamma$. There is a distinguished `blank' letter, $\beta$. 
\item It has a tape read/write head that starts on the square at position zero.

\item Finally, it has a transition function $\delta: (Q \setminus F) \times \Gamma \to Q \times \Gamma \times \{{\it L}, {\it R}\}$ that tells us what changes to make and which direction to go at each time step.
\end{itemize}
\smallskip

The machine takes as input a tape, which starts in some initial state~$T_0: \mathbb Z\to\Gamma$ with all but finitely many squares blank.

At every time step, if the head is in a final state, the machine halts.
Otherwise, let~$q \in Q \setminus F$ be the current state, and let $\gamma$ be the letter on the tape square that the head is currently on.
Let $\delta(q, \gamma) = (q', \gamma', d)$.
The machine changes to the new state~$q'$, writes the letter $\gamma'$ on the tape at the current position, and moves one square left or right depending on the direction $d$.
Once it's on the new square, it looks at the new state and letter, moves again, and so on.

If a tape square has the letter $\gamma$ on it, we will say that that square is in a {\sl tape state} $t(\gamma)$ if the head isn't on the square, and a {\sl head state} $h(q, \gamma)$ if the head is on the square in state $q$. If $q \in F$, we say that the square is in a {\sl final head state}.

The state of the square at position $i$ at a new time step is determined by the states of the three squares~$i-1$, $i$, $i+1$ at the previous time step, as long as none of those states are a final head state. The rules are in Figure~\ref{fusemove}.
\begin{figure}[h]
\begin{center}
\includegraphics[width=\hsize]{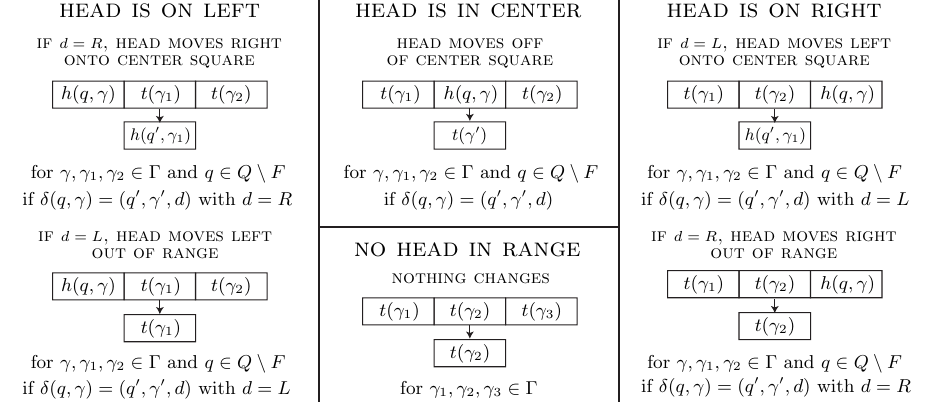}
\end{center}
\caption{How to find the new state of a tape square by looking at the state of it and its neighbours at the previous time step. There's only one tape head, so
these are the only four possibilities.}\label{fusemove}
\end{figure}

\end{subsection}

\begin{subsection}{Implementation as a sandpile via fuse circuits}
We build a periodic sandpile that calculates the state of each tape square in the Turing machine at each time step.
We start by cutting $\mathbb Z^3$ up into cubes of some large but unspecified side length $L$. We will label each cube by its position $(x, y, t)$: \begin{align*}c(x,y,t) = \{v \in \mathbb Z^3: {}&xL \le v_x < (x+1)L,\\&yL \le v_y < (y+1)L,\\
&tL \le v_z \le (t+1)L\}.\end{align*}
We're calling the third dimension $t$ because it keeps track of the flow of time in the simulated Turing machine, but it is really just the third space dimension in $\mathbb Z^3$ and has no relation to time in the sandpile. The cube at position $(x, 0, t)$ is going to calculate the state of the Turing machine square at position $x$ at time step $t$.

The sandpile configuration must be periodic in three dimensions by definition, but we only need to keep track of each tape square $x$ at each time step $t$. Therefore we only add connections between cubes with the same $y$-value, and we'll add extra chips only to the cubes on the $y = 0$ layer. Every other layer will be stable forever and won't affect the calculation on the $y = 0$ layer.\footnote{In fact, we don't need the third dimension at all except to allow wires to cross over, so we could build the same structure in a slab of height 3, $\{(x,y,z):0 \le y \le 2\}$.}

From now on we suppress $y$ in the notation, writing cubes as~$c(x, t) = c(x,y,t)$.

\begin{subsubsection}{Recording the cube state}If $S$ is a finite set of states, then we can assign a distinct bit string $s_1 \cdots s_b$ to every element $s \in S$, where $b = \lceil \log_2(\#S) \rceil$.

Each cube is going to store its current state, either $t(\gamma)$ for $\gamma \in \Gamma$ or $h(q, \gamma)$ for $q \in Q$, $\gamma \in \Gamma$,
so $S = \Gamma \cup (Q \times \Gamma)$ and we need $\lceil \log_2 \#\Gamma(1 + \#Q)\rceil$ bits.

We want to store one of these bit strings in each cube.
To do this, we put~$b$ pairs of state wires in every square of the lattice. We will call them $w_m(x, t)$ and $\overline{w}_m(x, t)$ for $m = 1, \ldots, b$.
At first, neither wire in a pair is toppled.
If~$w_m$ is toppled, it means that bit $m$ of the state is $1$; if the other wire $\overline{w}_m$ is toppled then bit $m$ is $0$.

If neither of the state wires have toppled yet, we say that bit $m$ is undecided, and if there are any undecided bits, we say that the cube is in an undecided state. Otherwise it's in the state corresponding to the bit string encoded in the wires. 

Under ordinary operation, at most one wire in a pair will topple.
A wire can't untopple, so a cube can go from the undecided state to a certain state $s \in S$, but once it has decided on its state, it can't go back to undecided or change its state.


\end{subsubsection}

\begin{subsubsection}{{\sc and} gates calculate the new state of the cube}\label{andgates}
Let $s_{-1}, s_0, s_{+1}$ be three states in $\Gamma \cup (Q \times \Gamma)$, at most one of which is a head state, and none of which is a head-final state $h(q,\gamma)$ for $q \in F$.
Let $\sigma = (s_{-1}, s_{0}, s_{+1})$.
We put a multi-way {\sc and} gate~$g_{\sigma}(x, t)$ for every cube $(x, t)$ and every tuple~$\sigma$ that satisfies the conditions. 

The inputs of the gate are connected to the state wires of neighbouring squares at the previous time step as follows.

For $k = -1, 0, 1$:
\begin{itemize}\item[]write $s_k$ as a bit vector $s_{k1}, \ldots, s_{kb}$, and then, for $m = 1, \ldots, b$:
\begin{itemize}\item[]if $s_{km} = 1$, connect an input to $w_m(x+k,t-1)$;\item[]otherwise, connect an input to $\overline{w}_m(x+k,t-1)$.\end{itemize}\end{itemize}
The {\sc and} gate fires if and only if the state wires in $c(x-1,t-1), c(x,t-1), c(x+1,t-1)$ match the bit patterns corresponding to $s_{-1},s_0,s_{1}$.

The output wire of the gate is connected to the state wires of the current square as follows. If the three squares at $x-1$, $x$, $x+1$ are in states $s_{-1}, s_0, s_{+1}$ at the previous time, then the square at position $x$ will enter some new state according to the transition rules in Figure~\ref{fusemove}. Call this new state $s'$. Let $s'_1, \ldots, s'_b$ be the bit string corresponding to~$s'$.
Then, for $m = 1, \ldots, b$:
\begin{itemize}\item[]if $s'_m = 1$, connect the output wire to $w_m(x, t)$;\item[]otherwise connect it to $\overline{w}_m(x, t)$.\end{itemize}
\end{subsubsection}

\begin{subsubsection}{If the cubes at time zero are set up correctly, then this fuse circuit simulates a Turing machine, until it halts}\label{ifzero}

Let $t > 0$. Suppose that the state wires of the three cubes $c(x-1,t-1)$, $c(x,t-1)$, $c(x+1,t-1)$ in the previous layer have been set to bit patterns that correctly represent the states of the squares of the Turing machine at positions $x-1$, $x$, $x+1$ at time step $t-1$.

If none of those is a final head state, then exactly one {\sc and} gate $g_\sigma$ will be able to activate in the cube $c(x, t)$, the one that's connected to the toppled wire in each pair of state wires in those three cubes. Every other {\sc and} gate in the cube is connected to at least one untoppled wire in some pair, so no other gate can activate.

When the {\sc and} gate fires it will change the state of $c(x,t)$ from undecided to $s'$, where $s'$ is the correct state of square $x$ at time $t$. So if the three cubes $c(x-1,t-1)$, $c(x,t-1)$, $c(x+1,t-1)$ in the previous layer eventually enter the correct state, then $c(x,t)$ will also eventually enter the correct state.

By induction, every cube $c(x,t)$ for $t > 0$ will eventually enter the state of square~$x$ at time step $t$, unless a final head state appears somewhere. If that happens, then every cube that sees that state at the previous time will remain undecided forever, because there are no {\sc and} gates that match the bit pattern of a final head state.

This means that if a cube enters a final head state, then every cube in its forward light cone~$c(x', t')$ for $|x-x'| \le t' -t$ will remain undecided, as in Figure~\ref{undecidedwave}.
If no cube enters a final head state, then the simulation will proceed correctly forever.
\begin{figure}[h]\includegraphics[width=\hsize]{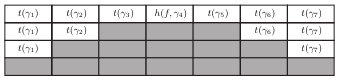}\caption{A wave of indecision for vertex halting. If a cube enters a final head state, every cube in its forward light cone stays undecided. Here $f \in F$ is a final head state and $\gamma_i$ are tape letters, and the gray cubes remain in an undecided state forever. In this case, the speed of light is $1$.}\label{undecidedwave}\end{figure}

\end{subsubsection}\begin{subsubsection}{Setting up the cubes at time zero}\label{initialization}
We saw in the last section that this fuse circuit will simulate a Turing machine, as long as the cubes are initialized correctly at time zero. But doing this presents a technical problem. We can set the state of finitely many cubes by adding non-periodic chips directly to their state wires, but we can't set all cubes $c(x, 0)$ for $x \in \mathbb Z$ with that approach, because we aren't allowed to add infinitely many non-periodic chips.

Fortunately, all but finitely many squares of the Turing machine's initial tape have the blank letter $\beta$ on them, so we will build a circuit that can initialize a half-line of cubes to the blank tape state $t(\beta)$ when we topple an input wire.

Hook up an {\sc or} gate $g_{\rm left}$ in every cube $c(x, t)$, with two input wires and $b+1$ output wires. One input wire is connected to an output wire of the same gate $g_{\rm left}$ in the cube $c(x+1,t)$ immediately to its right at the same time step, and another input wire is connected to nothing. Connect the remaining $b$ output wires of the {\sc or} gate up to the state wires in $c(x,t)$, as in the last section, 
 so that when the gate activates it sets the state of $c(x, t)$ to the blank tape state $t(\beta)$.

If this {\sc or} gate $g_{\rm left}$ activates in cube $c(x, t)$, then it sets the state of $c(x,t)$ to $t(\beta)$ and activates the $g_{\rm left}$ gate in cube $c(x-1, t)$. This sets off a chain of activations: by induction all cubes $c(x',t)$ for $x' < x$ will enter state $t(\beta)$.\footnote{Unless one of the other {\sc and} or {\sc or} gates tries to set the state to something else, in which case we will set both wires in a pair and the circuit will malfunction. We must take care to prevent that.}
We also set up a similar {\sc or} gate $g_{\rm right}$ where the chain of activations goes in the other direction.

We can set all but finitely many cubes $c(x, 0)$ to $t(\beta)$ just by putting two chips on the unconnected input wires of one $g_{\rm left}$ gate and one $g_{\rm right}$ gate, and then set the remaining cubes directly by toppling state wires.

To be precise, let $\emph{Non-blank} = \{0\} \cup \{x: S_0(x) \ne \beta\}$ be the set of cubes we want to set to something other than $t(\beta)$. Let $x_0 = \min \emph{Non-blank}$ and $x_1 = \max \emph{Non-blank}$. Topple all of the following $2 + b(x_1-x_0+1)$ wires:
\begin{itemize}
\item[] the unconnected wire of $g_{\rm left}$ in $c(x_0-1,0)$,
\item[] the unconnected wire of $g_{\rm right}$ in $c(x_0+1, 0)$,
\item[] and $b$ state wires in each cube $c(x',0)$ for $x_0 \le x' \le x_1$ so that:
\begin{itemize}
\item if $x' \ne 0$, the cube enters state $t(T_0(x'))$,
\item and $c(0,0)$ will enter state $h(q_0, T_0(0))$.
\end{itemize}
\end{itemize}
Recall that $q_0$ is the starting state of the Turing machine.

The ultimate result of that is that every cube $c(x',0)$ for $x' \ne 0$ will enter the tape state $t(T_0(x'))$, and $c(0, 0)$ will enter the head state $h(q_0, T_0(0))$.
That is the correct setup, so the cubes $c(x,t)$ for $t > 0$ will simulate the Turing machine.


\end{subsubsection}
\begin{subsubsection}{Note: we only need one non-periodic chip}
In the above scheme, we add many non-periodic chips.
We don't have to: we really only need to add \emph{one} non-periodic chip, the minimum number necessary to break the symmetry.

Here is how that works. Add an `initialization' wire to each cube $c(x,t)$. 
Connect that wire, via diodes, to every wire we want to topple initially.
In other words, connect it to $g_{\rm left}$ in $c(x+x_0-1,t)$ and $g_{\rm right}$ in $c(x+x_1-1,t)$, and connect it to $b$ state wires in each cube $c(x+x',t)$ for $x_0 \le x' \le x_1$ so that its state will become $t(T_0(x'))$ or $h(q_0,T_0(0))$ if $x' = 0$.

We connect each initialization wire to a finite number of wires in other cubes at a finite distance. That takes a finite amount of additional wiring per cube, but we can do that, although we might have to increase the side lengths of the cube to accommodate the additional wiring.

If we set up this wiring and then put one chip on the initialization wire in $c(0, 0)$, the Turing machine will be correctly initialized just as before, but this time with only one added non-periodic chip. We saw in section \ref{periodicsandpile} that a periodic sandpile is predictable, so we must add at least one non-periodic chip. Therefore this is the minimum possible.

However, the original setup has the advantage that we can change the initial tape just by changing the arrangement of non-periodic chips. So if the Turing machine we're simulating is universal, we can fix the periodic arrangement and take as input the finite list of non-periodic chips, and the resulting halting problems are still undecidable.
\end{subsubsection}
\end{subsection}

\end{section}

\begin{section}{Vertex prediction is undecidable}\label{vertexhalting}
The first problem in our list is not precisely a halting problem:

\smallskip
{\it P+f vertex prediction.}~Given a certain infinite periodic sandpile with finitely many chips added to it, does the vertex at the origin ever topple?

\begin{theorem} The p+f vertex prediction problem is undecidable.\label{pfvertexundecidable}
\end{theorem}

\begin{proof}
We will show that the halting problem for Turing machines can be reduced to vertex prediction: if we are able to solve this problem for any input, we can determine whether an arbitrary Turing machine halts.

In the previous section, we built a periodic+finite fuse circuit that simulates a Turing machine given a certain input.
All we need to do is add a little extra circuitry to make the origin topple when the Turing machine reaches a final state.

Put an wire in every cube $c(x, t)$, called the {\sl alarm wire}. Connect it to the alarm wires in adjacent cubes.
Put one {\sc and} gate in every cube for every $f \in F$ and~$\gamma \in \Gamma$, with inputs connected to the state cells so that the gate will activate if and only if the cube is in state $h(f, \gamma)$.
Connect the output to the alarm wire in that cube.
Finally, translate the whole sandpile so that one alarm wire goes through the origin.

Suppose the Turing machine being simulated enters a final state $f \in F$. It does so at some time step $t \ge 0$ with the head in some position $x \in \mathbb Z$.
The fuse circuit will simulate the Turing machine, and eventually the storage cells in the cube~$c(x, t)$ will be set to a bit pattern corresponding to some state~$h(f, \gamma)$. That will set off the alarm wire in that cube, which will set off the alarm wire in every cube, including the one that goes through the origin. So the origin topples.

If the Turing machine never enters a final state, then the alarm wires won't go off, and the origin will never topple.
So the vertex prediction problem for this sandpile~is equivalent to asking whether the Turing machine being simulated halts.

The input Turing machine is arbitrary, and the question of whether an arbitrary Turing machine halts is undecidable: there is no terminating algorithm that can solve the problem in all cases. So p+f vertex prediction is also undecidable.
\end{proof}

\end{section}

\begin{section}{Global halting: The lazy automaton}\label{lazy}
The second problem in our list is:

\smallskip
{\it P+f global halting prediction.} Given a periodic sandpile with finitely many chips added, are there infinitely many topplings in total?

\smallskip

To reduce Turing halting to global halting prediction, we need a sandpile circuit for which the total number of topplings $\sum_{v \in \mathbb Z^3} N(v)$ is finite if the Turing machine halts, and infinite if it does not.
The naive circuit in section~\ref{buildingturing} doesn't work for this: the initialization procedure in section~\ref{initialization} topples state wires in all cubes $c(x, 0)$, so we have already exceeded our budget just by initializing the tape.

Fortunately, we don't really need to initialize the whole tape at the very start. The head of a Turing machine can only move at speed one, so we can initialize only one blank square per time step in each direction
and still be safe, as long as we have
some way of stopping the initialization if the
machine halts.

We introduce the idea of a lazy cellular automaton and see how to implement it as a fuse circuit with the same halting behaviour. After that, we will build a lazy automaton which halts if and only if a Turing machine halts.

\begin{subsection}{Lazy cellular automata}
A {\sl lazy cellular automaton} with radius $r$ is a four-tuple $(S, \lambda, S_0, R)$, where $S$ is the set of states, $\lambda$ is an extra undecided or `lazy' state that is not in $S$, the initial condition of each square is $S_0: \mathbb Z \to S \cup \{\lambda\}$, and $R$ is the set of {\sl partial rules}, to be explained later.

The lazy automaton has an initial state and evolves with time. Unlike the sandpile or fuse circuit model, it is a synchronous model, where the outcome of the system depends crucially on the timing of events.

Let the state at time~$t$ be denoted by $S_t: \mathbb Z \to S \cup \{\lambda\}$. For a regular cellular automaton, the state of a square at time $t$ is uniquely determined by the states of neighbouring squares at time $t - 1$. This is not true for a lazy cellular automaton; instead it is equipped with a set of {\sl partial rules} $R$, which may not cover every case, and may even be contradictory.

A partial rule is a pair $\pi \to s$, where $s$ is a state in $S$ and  $\pi$ is a {\sl pattern}: a sequence $(\pi_{-r}, \ldots, \pi_{r}) \in (S \cup \{*\})^{2r+1}$.
Here $*$ is a special wildcard symbol which is not in $S$ and can appear only in patterns.
A pattern $\pi = (\pi_{-r}, \ldots, \pi_r)$ {\sl matches} the states of the automaton at position $x$ and time step $t$ if,
for all $k = -r, \ldots, r$:
\begin{itemize}
\item[]
either~$\pi_k = *$ is the wildcard symbol,
\item[] or $S_t(x) = \pi_k$. \end{itemize}

The automaton evolves as follows. For $x \in \mathbb Z$, $t \ge 0,$ let $R_{x, t}$ be the set of rules that match at position $x$ and time step $t$.
Then the new state of vertex $x$ at $t+1$ is $$\Sigma_{t+1}(x) = \begin{cases}s&\text{if there is exactly one rule $\pi \to s$ in $R_{x,t}$, or}\\\lambda&\text{if the set $R_{x, t}$ is empty.}\end{cases}$$

If there are no rules that match, that square enters a lazy state. If there are two or more rules with matching patterns at some time $t \ge 0$, that's not good: the whole automaton {\sl malfunctions} at that time, and the state is undefined for future times.

The automaton is said to {\sl halt} if it never malfunctions and all but finitely many squares are lazy at some time.
We will hook up a fuse circuit that works like this.

How do we simulate a synchronous automaton with an asynchronous sandpile? As in the Turing machine simulation from section~\ref{buildingturing}, we do it by repurposing one of the space dimensions as a time dimension, so that the state of the automaton at a certain time step $t$ is computed by the wires in a slab $\{(x, y, z): tn_z \le z < (t+1)n_z\} \subseteq \mathbb Z^3$.\end{subsection}

\begin{subsection}{Implementation as a fuse circuit}\label{fusecircuit}
We start with the implementation of the Turing machine back in section~\ref{buildingturing}, and modify that. As before, we start by cutting $\mathbb Z^3$ up into large cubes, which we label $c(x,t)=c(x,y,t)$. We assign a distinct bit string $s_1 \cdots s_b$ to every state $s \in S$ in the automaton, where $b = \lceil \log_2(\#S) \rceil$. The lazy state doesn't get a bit string.

We put $2b$ state wires in each cube as before. If the state wires are toppled in a bit pattern corresponding to a state in $S$, then we say the cube is in that state. If the state wires \emph{never} topple, so that the state remains undecided forever, we say that that cube is in the lazy state.

\begin{subsubsection}{{\sc and} gates implement the partial rules}
We put a multi-way {\sc and} gate~$g_{\rho}(x, t)$ for every cube $(x, t)$ and every partial rule~$\rho = (\pi_{-r}, \ldots, \pi_r) \to s$.

The inputs of the gate are connected to the state wires of neighbouring squares at the previous time in a way that depends on the pattern $\pi$. This is just a slightly more complicated version of section~\ref{andgates} where we don't always connect to every pair.

For $k = -r, \ldots, r$:
\begin{itemize}
\item If $\pi_k = *$, don't connect anything.
\item If $\pi_k \in S$, write $\pi_k$ as a bit vector $\pi_{k1} \cdots \pi_{km}$, and then:

For $m = 1, \ldots, b$,
 
\begin{itemize} \item[] if $\pi_{km} = 1$, connect an input to $w_m(x+k, t-1)$ \item[]otherwise, connect an input to $\overline{w}_m(x+k, t-1)$.\end{itemize}
\end{itemize}

Write $s = s_1 \cdots s_b$. For $m = 1, \ldots, b$:\begin{itemize}
\item[] connect an output to $w_m(x, t)$ if $s_m = 1$, or \item[] connect an output to $\overline{w}_m(x, t)$ if $s_m = 0$.
\end{itemize}

The {\sc and} gate $g_\rho(x, t)$ activates if and only if~$\pi$ matches at time $t-1$ at $x$.
When it fires, it sets one bit in each of the pairs of state wires so that the cube enters state~$s$. If no patterns match, then no gates activate, and the cube will be in the lazy state.\end{subsubsection}

\begin{subsubsection}{Relationship to the automaton}
We claim that this fuse circuit works just like a lazy cellular automaton,
 unless the automaton malfunctions.

Here is how the correspondence works. In a sandpile, and therefore in a fuse circuit, there is no concept of synchronization, but we can distinguish between wires that eventually topple and wires that never do.

If, in each pair $w_m(x, t), \overline{w}_m(x, t)$ in a cube $(x, t)$, exactly one of the wires topples, according to the bit pattern $s_1 \cdots s_b$, then the corresponding square $x$ of the cellular automaton will enter the state corresponding to that bit pattern at time $t$. If any of the pairs are lazy forever, then the square will be in the lazy state.

The {\sc and} gate $g_\rho(x, t)$ activates if and only if cubes $c(x-r, t-1), \ldots, c(x+r,t-1)$ eventually get into states that match the partial rule $\rho$. If they do, the {\sc and} gate will fire and the current cube will be set to the appropriate bit string. If none of the patterns match, none of the {\sc and} gates will ever fire, and none of the state bits will ever be set for that cube, so it will be lazy forever.

This is exactly the behavior of the lazy automaton, as long as at most one {\sc and} gate fires in every cube.
%
If two or more gates in the same cube fire, it's likely that both wires~$w_m, \overline{w}_m$ in at least one pair will be toppled, but this cannot happen unless more than one partial rule matches, in which case the automaton malfunctions. \end{subsubsection}

\begin{subsubsection}{If the automaton halts, only finitely many {\sc and} gates fire}
Each matching partial rule activates an {\sc and} gate, which sets the state of one cube and topples wires in that cube and its neighbours, but nothing more distant. If the lazy automaton halts, finitely many wires will topple in total, so the total number of topplings will be finite and the sandpile will halt. If the automaton doesn't halt, the state will be set in infinitely many cubes, so infinitely many {\sc and} gates will fire and the sandpile will not halt.
Let's write all this down as a lemma:

\begin{lemma}If there is no malfunction, the lazy automaton halts if and only if finitely many {\sc and} gates in the corresponding fuse circuit activate.\end{lemma}

When we implement the fuse circuit as a sandpile, the sandpile halts if and only if finitely many gates in the circuit ever fire. Therefore, we have the corollary:

\begin{corollary} {\it If there is no malfunction, the lazy automaton halts (in the sense that all states are $\lambda$ at some finite time) if and only if the corresponding sandpile halts (in the sense that the sandpile eventually stabilizes).}
\end{corollary}

The next step is to come up with a lazy automaton that simulates a Turing machine and that halts if and only if the Turing machine halts.
\end{subsubsection}
\end{subsection}
\end{section}

\begin{section}{A lazy automaton that simulates a Turing machine}\label{evenlazier}
\def\lw{{\blacktriangleleft}}
\def\rw{{\blacktriangleright}}

We are going to construct a lazy automaton that simulates a Turing machine. We must specify the set of states, the radius of the neighbourhood, the rules, and the initial condition.

\begin{subsection}{The set of states}

Our set $S$ of states has one state $t(\gamma)$ for every tape letter $\gamma \in \Gamma$, and one state $h(q, \gamma)$ for $q \in Q$ and $\gamma \in \Gamma$, plus two other states:~$\lw$, the left wavefront symbol, and $\rw$, the right wavefront symbol.

We set the radius of the automaton to $2$. That is, the new state at position $x$ depends on the \emph{five} squares $x- 2, x-1, x, x+1, x + 2$ at the previous time step.

\end{subsection}
\begin{subsection}{Partial rules}
The partial rules of the lazy automaton are summarized in the following figure, Figure~\ref{lazyrule}.
\begin{figure}[th]\label{lazyrule}
\begin{center}\includegraphics[width=\hsize]{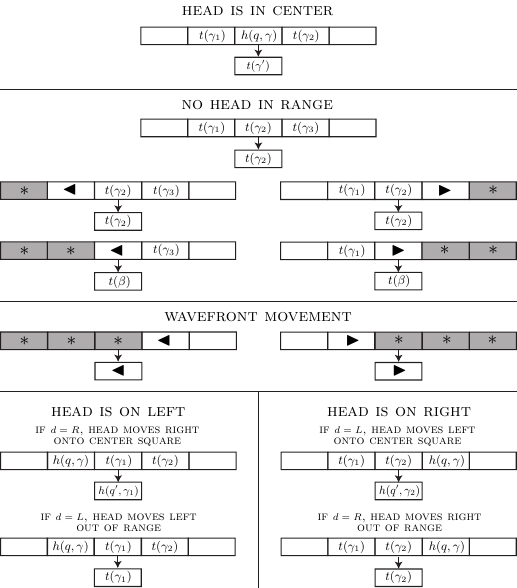}\end{center}\caption{The families of rules for the lazy automaton that simulates a Turing machine. Here $\gamma, \gamma_1, \gamma_2, \gamma_3$ are tape letters, $q \in Q \setminus F$ is a head state which must be non-final, and if $q$, $\gamma$ are present in a rule family, then $\delta(q, \gamma) = (q', \gamma', d)$. Gray squares are wildcards. Empty white squares are allowed to be any state other than a final head state.}\end{figure}

Each part of the figure represents a family of rules. Gray squares are wildcards. White squares can be any state in $S$ other than a final head state.
%
%
For example, the top family has $\#\Gamma^3(\#Q - \#F)(\#\Gamma(1+\#Q - \#F) + 2)^2$ rules, because we choose three tape letters $\gamma, \gamma_1, \gamma_2$, one non-final state $q \in Q \setminus F$, and two states that aren't final head states. Compare to Figure~\ref{fusemove}.

Because we have set the radius to~$2$, a partial rule is a pair $\pi \to s$ for $\pi \in (S \cup \{*\})^{5}$ and $s \in S$. We will write these as
$(a, b, \underline{c}, d, e) \to f$,
where the middle term is underlined to remind us which square is getting the new state.

We describe the families in detail in the rest of this section.
Let $\tau \subsetneq S$ be the set of states that aren't final head states, $\tau = \Gamma \cup ((Q \setminus F) \times \Gamma) \cup \{\lw, \rw\}$. 
This set has cardinality $\#\Gamma + \#\Gamma(\#Q-\#F) + 2$.
We will start with the rules that don't depend on the direction of motion of the Turing machine.

\begin{subsubsection}
{Rules when the head is in the center}\

\smallskip
{\bf If the head is on the center square, it changes the tape letter to $\gamma'$ and moves off of that square.}
When the head is on the center square, it's going to leave that square at the next time step, and it's going to change the tape letter to $\gamma'$ where $\delta(q, \gamma) = (q', \gamma', d)$. So for all $a, b \in \tau$, $\gamma, \gamma_1, \gamma_2 \in \Gamma$, $q \in Q \setminus F$, we add the rule:
$$(a,t(\gamma_1),\underline{h(q, \gamma)},t(\gamma_2),b) \to t(\gamma').$$

\end{subsubsection}
\begin{subsubsection}{Rules when no head is in range}\ 

\smallskip
All five families of rules in the second group are variations on a basic idea:

\smallskip
{\bf If all the squares within distance one are tape squares, the tape letter stays the same.}
For $a, b \in \tau$, $\gamma_1, \gamma_2, \gamma_3 \in \Gamma$,
$$(a,t(\gamma_1),\underline{t(\gamma_2)},t(\gamma_3),b) \to t(\gamma_2).$$

\smallskip
So far so good.
Now, if we can see a left wavefront symbol, we'd like to pretend that everything on the left side of it is a blank tape square. That gives us two more families of rules.

\smallskip
{\bf If we see a left wavefront to our left, we treat it as if the wavefront, and everything to its left, is a blank square.}
For $b \in \tau$, $\gamma_2, \gamma_3 \in \Gamma$,
\begin{align*}
(*,\lw,\underline{t(\gamma_2)},t(\gamma_3),b)
&\to t(\gamma_2)\\
(*, *, \underline{\lw},t(\gamma_3),b)
&\to t(\beta)
\end{align*}

\smallskip
We want to do the same thing for the right wavefront:

\smallskip
{\bf If we see a right wavefront to our right, we treat it as if the wavefront and everything to its right is a blank square.} For $a \in \tau$, $\gamma_1, \gamma_2 \in \Gamma$,
\begin{align*}
(a,t(\gamma_1),\underline{t(\gamma_2)},\rw,*)
&\to t(\gamma_2)\\
(a,t(\gamma_1),\underline{\rw},*,*)
&\to t(\beta)
\end{align*}
\end{subsubsection}

\begin{subsubsection}{Wavefront movement}
This group consists of two families of rules, one for the left wavefront and one for the right wavefront.

We want to move the left wavefront one square left each time step, provided that the square to its right is neither lazy nor in a final head state. If either of those two things is true, then our simulated Turing machine has halted and we are trying to stop the lazy automaton. So, we add the following family of partial rules:

\smallskip
{\bf If the left wavefront is one square to the right, and the square behind it is neither lazy nor in a final head state, then change the state to $\lw$.}

For $b \in \tau$,
$$(*, *, \underline{*}, \lw, b) \to \lw.$$

\smallskip
And the same thing for the right wavefront:

\smallskip
{\bf If the right wavefront is one square to the left, and the square behind it is neither lazy nor in a final head state, then change the state to $\rw$.}

For $a \in \tau$,
$$(a,\rw,\underline{*},*,*)\to\rw.$$

\end{subsubsection}
\begin{subsubsection}{Rules that depend on the direction of motion} If we see a head state at distance~$1$, we must decide whether the head moves onto our square or goes in the other direction and leaves us alone, by using the transition function $\delta$.

These rules are the same as the ones for the non-lazy Turing machine back in Figure~\ref{fusemove}, except that we require the squares at distance two to be in a valid state.

\smallskip {\bf If there's a head and it's moving left, then, if it's to our right, enter the correct head state; if it's to our left, stay in our current tape state.}

For $a, b \in \tau$, $\gamma, \gamma_1, \gamma_2 \in \Gamma$, $q \in Q \setminus F$, if $\delta(q,\gamma) = (q', \gamma', d) \text{ with }d = L$, then
\begin{align*}
(a,t(\gamma_1),\underline{t(\gamma_2)},h(q,\gamma),b)
&\to h(q', \gamma_2).\\
(a,h(q,\gamma),\underline{t(\gamma_1)},t(\gamma_2),b)
&\to t(\gamma_1)
\end{align*}

{\bf If there's a head, and it's moving right, then, if it's to our right, stay in our current tape state; if it's to our left, enter the correct head state.}

For $a, b \in \tau$, $\gamma, \gamma_1, \gamma_2 \in \Gamma$, $q \in Q \setminus F$, if $\delta(q,\gamma) = (q',\gamma',d)$ with $d = R$, then
\begin{align*}
(a,t(\gamma_1),\underline{t(\gamma_2)},h(q,\gamma),b)
&\to t(\gamma_2)\\
(a,h(q,\gamma),\underline{t(\gamma_1)},t(\gamma_2),b)
&\to h(q',\gamma_1)
\end{align*}

\end{subsubsection}

\begin{subsubsection}{Initial condition}Let $x_0$ be the first non-blank square in the initial condition of the Turing machine, or $-3$ if it would otherwise be larger: $x_0 = \min (\{-3\} \cup \{x: S_0(x) \ne \beta\})$. 
Let $x_1 = \max(\{+3\}\cup\{x: S_0(x) \ne \beta\})$.

We set the initial condition of the lazy automaton as in Figure~\ref{initialcondition}. All squares at positions $x \le x_0-2$ and $x \ge x_0+2$ start out in the lazy state.

\begin{figure}[h]\includegraphics[width=\hsize]{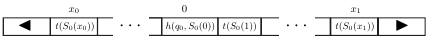}\caption{Initial condition for the lazy automaton.}\label{initialcondition}\end{figure}

The left wavefront symbol $\lw$ starts at position $x_0 \le -3$ and moves left, and the right wavefront symbol $\rw$ starts at position $x_1 \ge +3$ and moves right.

The Turing machine head moves at speed one in the automaton, so it will always be at least three cubes away from both the left and right wavefront, and
the distance between the two wavefronts is always at least seven, so no cube will ever see both wavefronts in its neighbourhood.
That simplifies the rules quite a bit.
\end{subsubsection}

\begin{subsubsection}{Wave of indecision}Everything in this system happens at speed one, until a cube in a lazy or final state appears between the two wavefronts. If it does, it will generate a wave of indecision as in Figure~\ref{undecidedwave}, but at speed \emph{two}.

Here is an example. Let $(Q,q_0,F,\Gamma,\beta,\delta)$ be a Turing machine with two states, $Q = \{q_0, q_1\}$. The second state is final, $F = \{q_1\}$. Let the alphabet be $\Gamma = \{\beta, 1\}$. Let $\delta(q_0, 1) = (q_0, \beta, R)$ and $\delta(q_0, \beta) = (q_1, \beta, L)$, so that the machine travels right, erasing a string of $1$s, until it hits a blank square, moves one square left, and halts. 
Let the initial tape be $\lw, \beta, \beta, 1, 1, 1, \rw$.
Then the lazy automaton evolves as in Figure~\ref{lazyexecute}. The horizontal axis is space, and the vertical axis is time.

\begin{figure}[h]\includegraphics[width=\hsize]{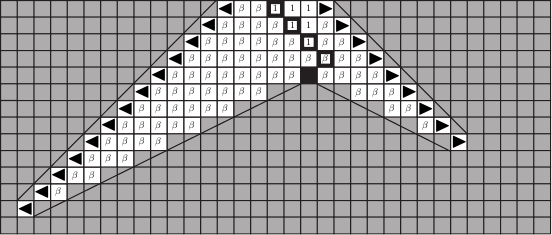}\caption{A simulation of a Turing machine by a lazy automaton. Gray squares are lazy, thin outlines are tape states, and heavy outlines are head states $h(q_0, -)$. The solid black square is the final head state $h(q_1, \beta)$. The lines show the speed of the wavefronts and the wave of indecision.}\label{lazyexecute}\end{figure}

The Turing machine halts at time 4, but the lazy automaton doesn't halt until the wave of indecision catches up with the left wavefront at time 13. This is true in general: if the Turing machine halts at position $x$ and time step $t$, it will take another~$t\pm x+3$ time steps for it to get rid of the wavefronts.

\end{subsubsection}

\end{subsection}
\begin{subsection}{This is a Turing machine simulation}

We can now prove that the global halting problem is undecidable by simulating this lazy automaton with a fuse circuit built in the sandpile.

\begin{theorem}The p+f global halting problem is undecidable.\label{pfglobalundecidable}
\end{theorem}

\begin{proof}
Given a Turing machine and the input for it, we construct a lazy automaton as above.
The automaton simulates a Turing machine until it enters a final state.
There are no partial rules for any final head states, so at the next time step every square within distance two becomes lazy.

A lazy state that appears between the two wavefronts will create a speed-two wave of indecision as above. It will eventually catch up to the speed-one wavefronts and halt the automaton.
On the other hand, if the Turing machine never enters a final state, no final or lazy state will appear between the wavefronts, so they will travel to the left and right forever, and the lazy automaton will not halt.

Then we build a periodic fuse circuit that simulates that lazy automaton as in section~\ref{fusecircuit}, and finally we build a periodic sandpile that operates in the same way as the fuse circuit, as in sections~\ref{oneshotsection} and~\ref{buildingturing}.

The lazy automaton halts if and only if the~Turing machine it is simulating halts, and the fuse circuit halts globally if and only if its automaton halts, so the global halting problem for periodic+finite sandpiles is undecidable.
\end{proof}

\smallskip

Note: Figure~\ref{lazyexecute} shows the evolution of the simulated lazy automaton with time, but it also shows the \emph{final} states of the sandpile cubes $c(x,0,t)$ in three-dimensional \emph{space}, where $x$ increases to the right and $z$ (aka ``$t$'') increases downward.

Time in the lazy automaton is $z$ in the sandpile simulation. Each time step corresponds to a line of cubes $c(x,0,t), x \in \mathbb Z$ in the sandpile.
To see what happens in the next time step of the automaton, we must move one period in the direction of positive $z$ and look at the next line of cubes.

\end{subsection}
\end{section}

\begin{section}{Local halting is undecidable}\label{localhalt}
We continue to the third halting problem.

\smallskip
{\it P+f local halting prediction.} Given a periodic sandpile with finitely many chips added, does the origin topple infinitely many times?

\smallskip
We begin with the sandpile that we built in section \ref{buildingturing} to simulate a Turing machine. Every vertex topples at most once, so the operation of the Turing machine itself will not make the origin topple infinitely many times.
We will set up a sandpile configuration that will topple every vertex infinitely often, and that can be set off by a wire. There will be a large region which is far enough away from criticality that the sandpile circuits we built in section~\ref{howgates} can function normally inside it.

Suppose that we have a cubical periodic domain in $\mathbb Z^3$ with large unspecified period $n$. We break the vertices in the periodic domain into four types.
A vertex $(x, y, z)$ is a {\sl face vertex} if exactly one coordinate is zero modulo $n$, an {\sl edge vertex} if two are, and a {\sl corner vertex} if all three are.
All other vertices are {\sl body vertices}.

Construct the Turing machine cellular automaton as in section \ref{buildingturing}, with one cell per Turing machine square per time tick, in such a way that the sandpile circuit:
\begin{itemize}
\item adds at most two chips to any body vertex that isn't part of a wire or gate, 
\item adds at most one chip to any face vertex that isn't part of a wire or gate,
\item and it adds no chips to any edge or corner vertex.
\end{itemize}
We can do that by making sure that no part of the sandpile circuit is near any edge or corner vertex, that the only sandpile circuit elements that go near the face vertices are straight wires running perpendicular through the face, that the bends in the wires are not too close together, and that there are a few squares of space between the gates.
None of the gates add more than two chips to a vertex that isn't part of the circuit.

Second, we put chips on every vertex without chips on it already: five on every edge and corner vertex, four on every face vertex, and three on every body vertex.  After that, all the edge and corner vertices will have five chips, every face vertex will have at least four, and every body vertex will have at least three.

Finally, add alarm wires and {\sc and} gates as in Section~\ref{vertexhalting}, but instead of connecting all the alarm wires together, connect each wire to the unique corner vertex in its cell.

The Turing machine operates as in section \ref{buildingturing}. As before, if the Turing machine reaches a final state, then one of the alarm wires will topple and set off a corner vertex. Once that happens, every single vertex in $\mathbb Z^3$ will topple, by the following lemma:

\begin{lemma} \label{bomb} If a sandpile on $\mathbb Z^3$ has at least:
\begin{itemize}
\item five chips on every vertex $(x, y, z)$ with at least two coordinates zero,
\item four chips on every vertex with one coordinate zero,
\item and three chips on every vertex,
\end{itemize}
and the origin topples at least once, then every vertex topples at least once.
\end{lemma}

\begin{proof} We will prove that all vertices $(x, y, z)$ with $x, y, z \ge 0$ topple and then extend this by symmetry to all eight octants.

We proceed by induction on the sum $x+y+z$ of the coordinates.
The origin topples, so that covers the base case $x+y+z = 0$.
Suppose $n \ge 1$, and all the vertices~$(x, y, z)$ with $x, y, z \ge 0$ and $x+y+z < n$ topple.

Let~$(x, y, z)$ be a vertex in $x, y, z \ge 0$ with $x+y+z = n$. There are three possibilities:

\smallskip
{\it The vertex has only one positive coordinate, say $x$.} Then it starts with at least five chips. By the induction hypothesis, it has at least one neighbour that topples, namely $(x-1, y, z)$. That neighbour adds one chip to it when it topples, so at that time it will have at least six chips on it, and it has to topple eventually by completeness.

\smallskip
{\it The vertex has two positive coordinates, say $x, y$.} In that case, the vertex starts with at least four chips. It has at least two neighbours that topple, $(x-1,y,z)$ and $(x,y-1,z)$. Once those two neighbours have toppled, the vertex will have at least six chips on it, and has to topple by completeness.

\smallskip
{\it The vertex has three positive coordinates.} Then it starts with at least three chips, and has at least three neighbours that topple: $(x-1,y,z)$, $(x,y-1,z)$, $(x,y,z-1)$. Again, once all three neighbours topple, the vertex has six chips and topples.

\smallskip
In all three cases, the vertex topples.
By induction, every vertex in the octant $\{(x, y, z): x, y, z \ge 0\}$ topples. The same argument works in every other octant.
\end{proof}


This and Lemma~\ref{oneinfinite} in the appendix give us the third theorem:
\begin{theorem}The p+f local halting problem is undecidable.\end{theorem}
\begin{proof}If we are given a Turing machine and an initial tape for it, we can build a periodic+finite sandpile with five chips on every edge and corner vertex, at least four chips on every face vertex, and at least three chips on every body vertex, for which the origin topples if and only if the Turing machine halts with that tape.

If the origin doesn't topple, the answer to the local halting problem is no. If the origin does topple, then the assumptions of Lemma~\ref{bomb} hold, so every vertex topples at least once. By Lemma~\ref{oneinfinite}, this means that every vertex topples infinitely often.

So the origin topples infinitely often if and only if the Turing machine halts at some time. We have reduced the halting problem for Turing machines to the p+f local halting problem, and so that problem is undecidable.
\end{proof}

\emph{Open problem.} A Turing machine can directly simulate the sandpile, and halt if the origin topples or if the sandpile globally halts, so p+f vertex and global halting problems are reducible to the halting problem for Turing machines.

Is p+f local halting reducible to the halting problem for Turing machines?
A priori, we only know that it's reducible to $\Pi^0_2$ in the arithmetical hierarchy, because the ability to evaluate formulas in $\Pi^0_2$ is equivalent to the ability to decide whether a Turing machine outputs infinitely many $1$s. The reduction is that we simulate the sandpile directly and output $1$ every time the origin topples. If we output infinitely many ones, then the sandpile does not halt locally; otherwise it does.

We could also ask the opposite question: if local halting is not reducible to $\Sigma_1$, then is local halting complete for $\Pi_2^0$, so that any problem in $\Pi_2^0$ can be reduced to the local halting problem for some periodic plus finite sandpile? This seems unlikely.

\end{section}

\appendix

\begin{section}{Elementary properties of sandpiles}
\begin{subsection}{The outcome of the sandpile process does not depend on the choices we make}\label{abelianproperty}

Recall from section \ref{abelianprop} that $G = (V, E)$ is a finite or countable graph, and $S_0: V \to \mathbb Z$ is the {\sl starting configuration} that tells us how many chips each vertex starts with. A {\sl toppling sequence} is a finite or infinite sequence of vertices in $V$, and a vertex $v \in V$ is {\sl unstable} if it has at least $\deg^+(v)$ chips on it.

A toppling sequence is said to be {\sl legal} if every vertex that's toppled is unstable, and it's {\sl complete} if every vertex that becomes unstable is eventually toppled later.


\begin{lemma}Given any finite or countable graph $(V, E)$ and any sandpile configuration on it, there exists a legal complete toppling sequence for that configuration.\label{legalcomplete}\end{lemma}

\begin{proof}
If $V$ is countable, let the vertices be numbered $1, 2, \dots.$ Look at the vertices in the infinite sequence $1, 2, 1, 3, 2, 1, 4, 3, 2, 1,$ and so on, and topple any vertex which is unstable when it is seen.
If $|V| = n$, then let the vertices be numbered $1, \ldots, n$, and look at them in the equally infinite sequence $1, 2, \ldots, n, 1, 2, \ldots, n, \ldots.$

In both cases, we never topple a stable vertex, and if a vertex is ever unstable, we will eventually check it again and topple it, so the result is legal and complete.
\end{proof}

We will think of a toppling sequence $a$ as a map $a: \iota_a \to V$, where~$\iota_a$ is some initial subset of $\{1, 2, \ldots\}$. Call $\iota_a$ the {\sl index set} of $a$.

If $a$ is a finite or infinite toppling sequence, let $N_a(v)$ be the number of times that the vertex $v \in V$ appears in the sequence, $N_a(v) = \#\{j \in \iota_a \mid a_j = v\}$.

If $a$ is legal and $b$ is complete, then the number of topplings at any vertex in the $a$ is less than or equal to the number of topplings at that vertex in $b$.
That is the content of the lemma below, due to Bj\"orner, Lov\'asz, and Shor \cite{bjornerlovaszshor}. The proof is adapted from Lemma 4.5 in Bond and Levine \cite{blev}.

\begin{lemma}\label{aplemma}
Fix a graph $(V, E)$ and a starting configuration $S_0: V \to \mathbb Z$. Then:

\begin{itemize}
\item[i)]
If $a$ is a legal toppling sequence for the configuration and $b$ is a complete toppling sequence for the configuration, then $N_a(v) \le N_b(v)$ for all  $v \in V$.

\item[ii)]
There is an \emph{odometer function} $N: V \to \mathbb Z_{\ge 0} \cup \{\infty\}$ so that every legal complete toppling sequence topples exactly $N(v)$ times at vertex $v$.
\end{itemize}
\end{lemma}

\begin{proof}

If the length of $a$ is zero, then $N_a(v) = 0$ for $v \in V$ and we are done. Otherwise, we are going to construct an injective map $\pi: \iota_a \to \iota_b$ which assigns an index $\pi(i)$ in $b$ for every index $i$ in $a$, and has the property that $a_i = b_{\pi(i)}$. The existence of such a map immediately implies that $N_a(v) \le N_b(v)$ for every vertex $v$.

%

We argue by induction. Suppose we have already chosen $\pi(1), \ldots, \pi(i-1)$ so that the map $\pi$ is injective on $\{1, \ldots, i-1\}$ and $a_{i'} = b_{\pi(i')}$ for $1 \le i' < i$. If the length of~$a$ is $i - 1$, then $\pi$ is already an injective map from $\iota_a$ to $\iota_b$ with the desired properties. Otherwise, let $\pi(i)$ be the smallest unused index $j$ in $\iota_b$ with $a_i = b_j$: $$\pi(i) = \min\{j \in \iota_b: b_j = a_i\ \text{and}\  j \notin \pi(\{1, \ldots, i-1\})\}.$$
We never reuse an index, so
this definition will give us the promised injective map, if there is always at least one unused index $j$ with $b_j = a_i$.
We will prove that there is.

Let $M = \max\{\pi(1), \ldots, \pi(i-1)\}$ be the maximum index used so far.

\medskip
\emph{If there is an unused index $j \le M$ with~$b_j = a_i$}, then we can extend the map by setting $\pi(i)$ to the first such index.

\medskip
\emph{If there is no index with $j \le M$ and $b_j = a_i$,} then we can use the completeness of $b$ to prove that there is an index $j$ with $b_j = a_i$ with $j>M$. Such an index is necessarily unused, because $j$ is larger than the maximum of every used index.

The proof goes as follows. Put chips on the graph so that it is in its starting configuration, and then topple at every vertex that has been used so far: at $b_{\pi(1)} = a_1$, then $b_{\pi(2)} = a_2$, and so on up to $b_{\pi(i-1)} = a_{i-1}$. By assumption, the sequence $a$ is legal, so the next vertex $a_i$ has at least $\smash{\deg^+(a_i)}$ chips after those topplings.

Now topple at the unused vertices $b_j$ for $j \in \{1, \ldots, M\} \setminus \pi(\{1, \ldots, i-1\})$. By assumption, none of those vertices are equal to $a_i$, so there are no additional topplings at the vertex $a_i$ and the number of chips on it has not decreased.

At this point, we have toppled every vertex $b_1, \ldots, b_M$, and the vertex $a_i$ has at least $\deg^+(a_i)$ chips on it, so it is unstable. The sequence $b$ is complete, so there must be an index $j > m$ with $b_j = a_i$.

Therefore, there is always at least one unused index $j$ with $b_j = a_i$, and we can construct the injective map.
By the reasoning in the first paragraph, $N_a(v) \le N_b(v)\ \forall v \in V$, and we have proven the first statement of the lemma.

For the second statement, pick any legal complete toppling sequence $a$, and let~$N(v) = N_a(v)$. If $b$ is any other legal complete toppling sequence, two uses of the first part of this lemma tell us that $N_b(v) = N(v)$ for every vertex $v \in V$.
\end{proof}

\end{subsection}

\begin{subsection}{The odometer function is the solution of a minimization problem}

There is another characterization of the odometer function $N$ from Lemma~\ref{aplemma}.

\medskip

\emph{Definition.} If $f: V \mapsto \mathbb Z$ is a function on the vertices of the graph, let the graph Laplacian $\nabla^2 f$ be the function
\[\nabla^2 f(v) = \sum_{w \sim v} f(w) - f(v).\]

\begin{lemma}If $(V, E)$ is a connected, undirected graph, and the degree of every vertex is finite, then the odometer function $N: V \to \mathbb Z_{\ge 0}$ for a starting configuration~$S_0$ is the solution of the following minimization problem, if the set is nonempty:
\begin{equation}N(v) = \min\left\{f(v) \mid f: V \to \mathbb Z_{\ge 0},\,\forall v \in V: S_0(v) + \nabla^2 f(v) < \deg(v)\right\}.\label{minimuexpression}\end{equation}
Moreover, $N$ is an element of the above set.
If the set is empty, then $N \equiv +\infty$.\label{variat}
\end{lemma}

\begin{proof}Let a toppling sequence $a$ that contains each vertex $v$ at most finitely many times be called \emph{locally finite}. The degree of each vertex is finite by assumption, so there are only finitely many topplings on $v$ and its neighbours, and the number of chips on $v$ will eventually be constant. The final number of chips will be
$$S_0(v) - N_a(v) \deg v + \sum_{w \sim v} N_a(w) = S_0(v) + \nabla^2  N_a(v).$$
If this is strictly less than $\deg(v)$ for every vertex $v \in V$, then every vertex is eventually stable; otherwise, some vertex stays unstable forever. So a locally finite sequence~$a$ is a complete sequence if and only if $S_0 + \nabla^2 N_a < \deg v$.

Conversely, if $f \ge 0$ is a finite integer-valued function with $S_0 + \nabla^2 f < \deg v$, then there is a locally finite toppling sequence $a_f$ which topples at $v$ a total of $f(v)$ times. By the reasoning above, every vertex is eventually stable, so $a_f$ is complete. By Lemma~\ref{aplemma}, the odometer function is bounded above by $N_{a_f} = f$.

If the set of functions $f$ in $(\ref{minimuexpression})$ is nonempty, then the odometer function is bounded above by a finite function $f$. Let $a$ be a legal complete toppling sequence, which must be locally finite because $N \le f < \infty$ everywhere. By the reasoning above, $N = N_a$ is itself a finite nonnegative function with $S_0 + \nabla^2 N < \deg v$.

On the other hand, if the set is empty, then $N$ is not in it. So $N$ must not be finite everywhere, and it follows that $N(v) = +\infty$ everywhere by the following Lemma \ref{dichotomy}.
These two statements together give us the conclusion. \end{proof}

\begin{lemma}If $(V, E)$ is a connected graph, and the out-degree of every vertex is finite, the odometer function is either finite at every vertex or infinite at every vertex.\label{dichotomy}\end{lemma}

\begin{proof}Suppose $v$ is a vertex with $N(v) = +\infty$, and $w$ is a neighbour. Choose a legal complete toppling sequence. The vertex~$v$ appears infinitely many times in the sequence, and $w$ receives one chip every time that $v$ topples.
If~$w$ topples only finitely many times, then at some step of the sequence it will topple for the last time, but~$w$ will receive infinitely many chips afterward and eventually become unstable since the degree of $w$ is finite. This is impossible by completeness. Therefore, $N(w) = +\infty$.

If $N(v) = +\infty$, then $N(w) = +\infty$ for every out-neighbour $w$. By connectedness of the graph, either $N(v) = +\infty$ for every vertex or $N(v) < \infty$ for every vertex.
\end{proof}
\end{subsection}
\begin{subsection}{If the odometer function of a connected, undirected graph is always at least one, then it is infinite}

\begin{lemma}If the odometer function for a sandpile configuration on an undirected graph is always at least $1$, then it is infinite everywhere.\label{oneinfinite}\end{lemma}

\begin{proof}We may assume without loss of generality that the graph is connected.
By Lemma~\ref{variat}, the odometer function $N$ is the minimum of all nonnegative functions $f \ge 0$ with $S_0(v) + \nabla^2 f(v) < \deg v\ \forall v$, or it is infinite if there are no such functions.

Suppose that there is such a function. We saw in the proof of the lemma that $N$ itself satisfies the discrete inequality $S_0(v) + \nabla^2 N(v) < \deg(v)$. Set $f(v) = N(v) - 1$. Then $f \ge 0$, and the discrete Laplacian of $1$ is zero, so $S_0(v) + \nabla^2 f(v) = S_0(v) + \nabla^2 N(v) < \deg(v)$ for every vertex $v$, but $f < N$, which is a contradiction.

Therefore there is no such function and $N \equiv +\infty$.
\end{proof}
\end{subsection}
\begin{subsection}{The odometer function of a periodic configuration is the same as the odometer function on the torus graph}
Suppose we have a periodic sandpile configuration, as in Section~\ref{periodicsandpile}, with three periods $n_x$, $n_y$, $n_z$. We compare it to the sandpile on the finite torus graph $T = \mathbb Z^3 / (n_x \mathbb Z \times n_y \mathbb Z \times n_z \mathbb Z)$.

\begin{lemma}If $S_0$ is periodic with periods $n_x, n_y, n_z$, and $s_0$ is the corresponding configuration on the torus $T := \mathbb Z^3 / (n_x \mathbb Z \times n_y \mathbb Z \times n_z \mathbb Z)$, then the odometer function~$N$ on the lattice is the periodic extension of the odometer function $N_T$ on the torus.\label{sameastorus}\end{lemma}
\begin{proof}
We will prove first that the odometer function $N$ is periodic. Let $e_x = (n_x, 0, 0)$ be the period vector in the $x$ direction. Consider the translated odometer function $N_x = N(v - e_x)$. It is nonnegative, and $\nabla^2 N_x(v) = \nabla^2 N(v - e_x)$.

We will use Lemma \ref{variat}. The initial configuration is periodic, so $S_0(v) = S_0(v - e_x)$. Therefore,
$S_0(v) + \nabla^2 N_x(v) = S_0(v - e_x) + \nabla^2 N(v - e_x) < 2d$, and of course $N_x \ge 0$, so, by the variational characterization, the function $N_x$ is an upper bound for $N$.
Similarly, $N_{-x} = N(v + e_x)$ is an upper bound for $N$.

But then $N(v) = N_{-x}(v - e_x) \ge N(v - e_x) = N_x$. Putting those two inequalities together, we get $N = N_x$. Therefore $N$ is periodic in the $x$ direction with period $n_x$. The same reasoning tells us that it is periodic in other directions.

Now that we have proven that the odometer function $N$ on the lattice is periodic, we can factor through the quotient map to get a function $\tilde N: T \to \mathbb Z^3$. Then $\tilde N$ is an upper bound on $N_T$ by Lemma \ref{variat}, and the periodic extension of $N_T$ is an upper bound for $N$, so both of the inequalities are equalities and we have the result.
%
%
%
\end{proof}
\end{subsection}
\begin{subsection}{Guarantees about halting}
We will use the variational characterization of the odometer function to prove a general lemma about halting, which is used in the next section, and also justifies the claim in Section~\ref{haltingmentionedsection} that a sandpile with finitely many chips must halt (although there are easier proofs for that special case).

\begin{lemma}There is a nonnegative, integer-valued function $f: \mathbb Z^d \to \mathbb Z_{\ge 0}$ with $\sum_{v \in \mathbb Z^d} f(v) < \infty$ and $\nabla^2 f(v) < 2d-B$ when $||v||_1 \le R$ and $<d+1$ when $||v||_1 > R$.\label{differenceq}\end{lemma}

We will postpone the proof for a little bit. With this lemma, we are able to solve the global halting problem for a large class of configurations.

\begin{theorem}\label{dminustwoalwayshalt}If $S_0$ is a sandpile configuration on $\mathbb Z^d$ with $S_0(v) \le d-1$ for all but finitely many vertices $v$, it will halt globally.\end{theorem}
\begin{proof}Let $B := \max_{v \in \mathbb Z^d} S_0(v)$. Let $R := \max\{||v||_1: v \in \mathbb V,\ S_0(v) \ge d\}$.
Let~$f$ be the function from Lemma~\ref{differenceq}.
%
Then $S_0(v) + \nabla^2 f(v) < 2d$ for every vertex $v \in V$, and~$f$ is nonnegative, so by Lemma~\ref{variat}, the odometer function is bounded above by~$f$. Therefore, $\sum_v N(v) \le \sum_v f(v) < \infty$ is finite and the sandpile halts globally.
\end{proof}

We will now prove the lemma. The proof is a bit long, but elementary.
\begin{proof}[Proof of Lemma~\ref{differenceq}]

There is a decreasing function $h$ from $\{-1, 0, 1, \ldots\}$ to $\mathbb Z_{\ge 0}$ which is zero for sufficiently large arguments and satisfies the inequality
\[h(i-1) - 2h(i) + h(i+1) \le \begin{cases}2 - {B+1 \over d}&\text{ if }i \le R\\1&\text{ otherwise}.\end{cases}\] for all nonnegative integers $i$. If $B\le 2d-1$, then $h \equiv 0$ solves the difference inequality. Otherwise, we can write a satisfactory function explicitly as follows:
\[h(i) = \begin{cases}a + xi(i+1)/2 &\text{ for }{-1} \le i < R\\
(s-i)(s+1-i)/2 &\text{ for }R \le i \le s\\
0 &\text{ for } i > s,\end{cases}\]
where $x = \lfloor {2 - {B+1 \over d}} \rfloor < 0$ and $s = (R+1) (1 - x)$, and $a = (R+1)^2x(x-1)/2$, which is an integer. We leave it to the reader to check that this has all the required properties: it is decreasing, nonnegative, integer-valued, and satisfies the inequality.

Fix $v \in \mathbb Z^d$. Set $f(v) := h(||v||_1)$.
Let $m := \#\{w \sim v: ||w||_1 < ||v||_1\}$. Then $m \le d$, and $h$ is decreasing, so $(d-m)(h(i-1) - h(i+1))$ is nonnegative. Therefore,
\[\begin{aligned}
	\nabla^2 f(v) &= mh(i-1) - 2dh(i) + (2d-m) h(i+1)
	\\&\le d \times (h(i-1) - 2h(i) + h(i+1)).
\end{aligned}\]
This is less than $2d - B$ if $i \le R$, and less than $d + 1$ if $i > R$, so we are done.
\end{proof}
\end{subsection}
\begin{subsection}{There are no crossover gates in $\mathbb Z^2$}

In this section we classify some types of sandpile configurations on $\mathbb Z^2$ as `four-wire gates,' and prove that there is a planarity restriction on their behavior: a `crossover gate' cannot be constructed.
The argument that follows is essentially from Gajardo and Goles~\cite{gajardo}.

Let the {\sl basic four-wire configuration} on $\mathbb Z^2$ be the sandpile configuration with three chips in the set $W := \{v \in \mathbb Z^2: v_x v_y = 0\}$, and no chips anywhere else.

\begin{definition}A sandpile configuration on $\mathbb Z^2$ is a {\sl four-wire gate} if it is equal to the basic four-wire configuration except at finitely many vertices.
\end{definition}


As in the main text, a wire is a connected set of vertices which all have $2d-1 = 3$ chips on them.
A four-wire gate has four half-infinite wires which stretch off forever to the north, east, south, and west. We will call these the {\sl cardinal wires.}

We can send signals to this gate by toppling its wires, and the gate can send signals out by toppling some wires itself. It's intuitively plausible that such a gate can only send signals `to infinity' by toppling the wires. The next lemma proves that.

\begin{lemma}Let $N(v)$ be the odometer function of a four-wire gate. Then, for all but finitely many vertices, $N(v) = 0$ for $v \notin W$, and $N(v) \le 1$ for $v \in W$. Also, $N$ is constant on each cardinal wire, except for finitely many vertices.\end{lemma}
\begin{proof}Let $S_0$ be the starting configuration of the gate. If $h(v) = 1$ for $v \in W$ and $0$ for $v \notin W$, then the new configuration $S_0' := S_0 + \nabla^2 h$ is less than or equal to~$1$ for all but finitely many vertices. By Lemma~\ref{dminustwoalwayshalt}, it globally halts. So the odometer function~$N'$ of the new configuration is zero for all but finitely many vertices.

Lemma~\ref{variat} says that the new odometer function $N'$ satisfies $S_0' + \nabla^2 N' < 4$, but we defined $S_0' = S_0 + \nabla^2 h$, so we also have $S_0 + \nabla^2 (N' + h) < 4$ and, again by Lemma~\ref{variat}, we must have $N \le N' + h$, which is equal to $h$ for all but finitely many vertices. That immediately implies the first statement in the result.

The vertices in each cardinal wire all start with $3$ chips, so if any of them topples, they all do. Fix a cardinal wire. Either $N(v)$ is identically zero for all vertices on the wire, or $N(v) \ge 1$ for all vertices in the wire. In the first case, we are done: $N$ is zero on the wire.
In the second case, we have proven in the previous paragraph that $N \le h \le 1$ for all but finitely many vertices, so we must have $N(v) = 1$ for all but finitely many vertices in the cardinal wire. That is the second statement. \end{proof}

If we add a chip to a cardinal wire, it topples. The above lemma tells us that the place where we add the extra chip does not affect the odometer function, as long as it is far enough away from the origin\footnote{Suppose that the four-wire gate topples the wire anyway. Let $v$ be far enough away that $S_0(v) = 1$, and $N = h$ for $v$ and all its neighbours. Then $v$ starts with three chips, topples once and has two neighbours that topple once, so the final number of chips is $1$. Adding an extra chip brings that to $2$, which is not enough to cause another toppling.
If the gate doesn't topple the wire, replace it by one that does, which only increases the odometer function, and follow the same reasoning. The vertex must topple once, and not twice, so the position of the chip doesn't change the odometer.}. We say that we have added a chip `at infinity.'

Suppose that we look at the gate from very far away, ignoring what happens near the origin.
If we put chips on some of the cardinal wires `at infinity,' then the gate may cause some of the other cardinal wires to topple, and that is the only effect we can see at long range, by the previous lemma.

We can characterize the gate in these terms. For example, we could say that a three-input {\sc and} gate is one that topples the south wire if and only if the north, east, and west wires are toppled,
and it never topples the other three wires. It is possible to build a three-input {\sc and} gate, and we have already seen how to do it.

Can we build a crossover gate that allows two wire signals to cross over each other without interfering, so that our circuits don't have to be planar? Such a gate would allow us to transfer the sandpile circuits to $\mathbb Z^2$ without any planarity limitations.

\begin{definition}
A four-wire gate is a {\sl crossover gate} if toppling the east wire causes the west wire to topple, but not the north wire, and on the other hand, toppling the south wire causes the north wire to topple, but not the west wire.\end{definition}

We will prove that there is no configuration that behaves like this by comparing odometer functions and using the planarity of $\mathbb Z^2$.

\begin{lemma}There's no way to construct a crossover gate.\label{crossover}\end{lemma}

\begin{proof}Suppose that there is a crossover gate configuration. By adding diodes far from the origin (for example, at distance at least~$s + 2$ in Lemma~\ref{dminustwoalwayshalt}), we can prevent the east and south wires from being toppled by the gate. 
Suppose we have done this.

Let the configuration with one chip added to the east (south) wire at infinity be called $S_e$ (or $S_s$) with odometer functions $N_e$ (or $N_s$) respectively.
The assumptions say that $N_e = 1$ for all but finitely many sites on the east and west wires, but $N_e = 0$ for all but finitely many sites on the north and south wires.

Let $E := \{v: N_{e}(v) > N_{s}(v)\}$ and $F := \{x: N_{s}(v) > N_{e}(v)\}$. These two sets are disjoint, and $E$ contains all but finitely many sites in the east and west cardinal wires, while $F$ contains all but finitely many sites in the north and south cardinal wires.
Either the east and west cardinal wires are disconnected in $E$, or the north and south cardinal wires are disconnected in $F$; this follows from planarity. We can flip the gate diagonally, so suppose without loss of generality that the east and west wires are disconnected in $E$.

Let $E'$ be the component of $E$ that contains the west cardinal wire. Set $$N'_{e}(v) := \begin{cases}N_{s}(v)&v \in E'\\N_{e}(v)&v \notin E'.\end{cases}$$
Then we claim that $S_e(v) + \nabla^2 N_e'(v) < 2d$ for every vertex $v \in \mathbb Z^d$.

If $v \in E'$, then we have $N_e'(v) = N_s(v)$, and all the neighbours $w \sim v$ are either in $E'$ or $F$, so $N_e'(w) \le N_s(w)$. Therefore, $S_e(v) + \nabla^2 N_e'(v) \le S_s(v) + \nabla^2 N_s(v) < 2d$.
Here we use the fact that $E'$ does not contain any site in the east cardinal wire, so $v$ does not have the east chip at infinity on it, and $S_e(v) \le S_s(v)$.

If $v \notin E'$, then we have $N_e'(v) = N_e(v)$, and the neighbours $w \sim v$ will either have $N_e'(w) = N_s(w) < N_e(w)$ if they are in $E'$, or $N_e'(w) = N_e(w)$ if they are not in~$E'$. Again, we will have $S_e(v) + \nabla^2 N_e(v) \le S_e(v) + \nabla^2 N_e(v) < 2d$.

So $S_e(v) + \nabla^2 N_e'(v) < 2d$ for every vertex $v \in \mathbb Z^d$. By Lemma~\ref{variat}, we must have $N_e \le N_e'$, but this is a contradiction since $N_e' = N_s < N_e$ on the nonempty set $E'$.
Therefore, our assumption was false, and there is no crossover gate.
\end{proof}

This does not totally rule out the possibility of nonplanar signaling in sandpiles on $\mathbb Z^2$, but it does tell us that there is no way to send a nonplanar signal with wires, or in any way that can be converted to a wire signal.

\end{subsection}

\end{section}

\bibliographystyle{acm}
\bibliography{sandpile_turing_machine-arxiv}

\section*{Acknowledgements}

Thanks to my advisor, Lionel Levine, who suggested that I put this idea in the form of a paper and made extensive suggestions. Matt Farrell and Swee Hong Chan kindly read my drafts and gave valuable comments. Also thanks to the referees, who pointed out Moore and Nilsson's bound and made the remark at the end of section~\ref{localhalt}.

\end{document}